\documentclass[]{informs1}              % for a regular run
\usepackage{natbib}
 \NatBibNumeric
 \def\bibfont{\small}%
 \bibpunct[, ]{[}{]}{,}{n}{}{,}%

\usepackage{mathrsfs}
\usepackage{amsmath}
\usepackage{amssymb}
\usepackage{bbm}
\usepackage{tikz}
\usepackage{setspace}
\usepackage{url}
\usetikzlibrary{calc, shapes, backgrounds}
\usetikzlibrary{decorations.pathreplacing}
\usetikzlibrary{patterns}
\usetikzlibrary{arrows}

\usepackage[colorlinks=true,breaklinks=true,bookmarks=true,urlcolor=blue,
     citecolor=blue,linkcolor=blue,bookmarksopen=false,draft=false]{hyperref}

\def\EMAIL#1{\href{mailto:#1}{#1}}% When hyperref is used, otherwise outcomment 
\TheoremsNumberedThrough     % Preferred (Theorem 1, Lemma 1, Theorem 2)
\EquationsNumberedThrough    % Default: (1), (2), ...
\begin{document}
%%%%%%%%%%%%%%%%

%\setcounter{page}{1}
%\DOI{10.1287/xxxx.0000.0000}%


%\RUNTITLE{}

% \TITLE{Bundling Information Goods of Decreasing Value}
\TITLE{2-Approximation Algorithms for Perishable Inventory Control When FIFO Is an Optimal Issuing Policy}

\ARTICLEAUTHORS{%
\AUTHOR{Can Zhang}
\AFF{H. Milton Stewart School of Industrial and Systems Engineering, Georgia Institute of Technology, Atlanta, Georgia, 30332, \EMAIL{czhang2012@gatech.edu}}
\AUTHOR{Turgay Ayer}
\AFF{H. Milton Stewart School of Industrial and Systems Engineering, Georgia Institute of Technology, Atlanta, Georgia, 30332, \EMAIL{ayer@isye.gatech.edu}}
\AUTHOR{Chelsea C. White III}
\AFF{H. Milton Stewart School of Industrial and Systems Engineering, Georgia Institute of Technology, Atlanta, Georgia, 30332, \EMAIL{cwhite@isye.gatech.edu}}
} % end of the block

\ABSTRACT{%
% We consider a periodic-review, fixed-lifetime perishable inventory control problem where demand is a general stochastic process. The optimal solution for this problem is computationally intractable due to the ``curse of dimensionality''. In this paper, we first present a \textit{marginal-cost dual-balancing policy} for our problem. We then prove that a myopic policy under the so-called marginal-cost accounting scheme provides a lower bound on the optimal ordering quantity. By truncating the balancing policy with this specific lower bound, we present a \textit{truncated-balancing policy}. We prove that when first-in-first-out (FIFO) is an optimal issuing policy, both of our proposed algorithms admit a worst-case performance guarantee of two, i.e. the expected total cost of our policy is at most twice that of an optimal ordering policy. We further present sufficient conditions that ensure the optimality of FIFO issuing policy. Finally, we conduct numerical analyses based on real data and show that both of our algorithms perform much better than the worst-case performance guarantee, and the truncated-balancing policy has a significant performance improvement over the balancing policy.
We consider a periodic-review, fixed-lifetime perishable inventory control problem where demand is a general stochastic process. The optimal solution for this problem is intractable due to ``curse of dimensionality''. In this paper, we first present a computationally efficient algorithm that we call the \textit{marginal-cost dual-balancing policy} for perishable inventory control problem.  We then prove that a myopic policy under the so-called marginal-cost accounting scheme provides a lower bound on the optimal ordering quantity. By combining the specific lower bound we derive and any upper bound on the optimal ordering quantity with the marginal-cost dual-balancing policy, we present a more general class of algorithms that we call the \textit{truncated-balancing policy}. We prove that when first-in-first-out (FIFO) is an optimal issuing policy, both of our proposed algorithms admit a worst-case performance guarantee of two, i.e. the expected total cost of our policy is at most twice that of an optimal ordering policy. We further present sufficient conditions that ensure the optimality of FIFO issuing policy. Finally, we conduct numerical analyses based on real data and show that both of our algorithms perform much better than the worst-case performance guarantee, and the truncated-balancing policy has a significant performance improvement over the balancing policy.
}%

%\KEYWORDS{deterministic inventory theory; infinite linear programming duality; 
%  existence of optimal policies; semi-Markov decision process; cyclic schedule}
%\MSCCLASS{Primary: 90B05; secondary: 90C40, 90C90}
%\ORMSCLASS{Primary: Inventory/production: deterministic multi-item;
%  secondary: dynamic programming/optimal control: deterministic 
%  semi-Markov; programming: infinite dimensional}
%\HISTORY{Received November 20, 2003; revised March 8, 2004, and March 26, 2004.}

\KEYWORDS{perishable inventory; nonstationary correlated demand; approximation algorithms; optimality of FIFO issuing policy}
%\MSCCLASS{}
%\ORMSCLASS{Primary: ; secondary: }

\maketitle
%%%%%%%%%%%%%%%%%%%%%%%%%%%%%%%%%%%%%%%%%%%%%%%%%%%%%%%%%%%%%%%%%%%%%%

% Samples of sectioning (and labeling) in MOOR.
% NOTE: (1) all section levels end with a period,
%       (2) capitalization is as shown (sentence style, not title style).
%
%\section{Introduction.}\label{intro} %%1.
%\subsection{Duality and the classical EOQ problem.}\label{class-EOQ} %% 1.1.
%\subsection{Outline.}\label{outline1} %% 1.2.
%\subsubsection{Cyclic schedules for the general deterministic SMDP.}
%  \label{cyclic-schedules} %% 1.2.1
%\section{Problem description.}\label{problemdescription} %% 2.

% Text of your paper here

\section{Introduction.}
\label{sec:intro}

Perishable products are very common in practice. Typical examples include medical products such as blood and certain pharmaceuticals, and food products such as refrigerated meat and many dairy products. Unlike nonperishable products that can wait in inventory until they are used to satisfy demand, perishable products must be used within a short period of time, and will become outdated otherwise. Outdating can result in a significant amount of wastes and financial losses. For example, the number of platelets outdated in 2011 in the U.S. was approximately 321,000 units, which accounted for 12.8\% of all processed units (\citet{us2011}). Similarly, the total annual unsaleable costs in the food, beverage, health and beauty industries in the U.S. were estimated as \$15 billion, and about 17\% of these costs (over 2.5 billion dollars) were caused by outdating (\citet{grocery2008joint}). These facts underline the critical need for efficient inventory management policies for perishable products.

Our study is specifically motivated by a platelet inventory control problem faced by a local acute-care hospital, Hospital Alpha (name blinded). In Hospital Alpha, the demand for platelets mainly comes from cardiac surgeries, which account for more than 85\% of its platelet transfusion. In this case, the uncertainty of demand stems from two sources: 1) the number of surgeries performed per day, and 2) the amount of platelets needed per surgery. Such a compound structure of demand is common for many blood products. As such, the compound Poisson distribution, where a random (Poisson) amount of patients arrive at every time period and each patient consumes a random amount of blood products, has been widely assumed for modeling demand in the blood supply chain literature (e.g., \citet{gregor1982evaluation, kopach2003models, katsaliaki2008cost}). However, while simply assuming random arrivals is reasonable for some cases such as trauma patients, forecast information on the number of arrivals is often available for many other cases, especially for scheduled operations such as cardiac surgeries. In particular, most of those surgeries are scheduled days or even weeks in advance, thus the number of surgeries scheduled for each day is gradually revealed as time approaches. Although the compound structure of demand is widely considered in the blood inventory management literature, to our knowledge, the dynamically evolving forecast information on the number of arrivals is not formally captured.

Motivated by the platelet inventory control problem with evolving forecast information, in this paper, we study a periodic-review, fixed-lifetime perishable inventory control problem under a general demand process which can be nonstationary, correlated, and dynamically evolving over time. {Similar to many other perishable inventory focused studies, we consider first-in-first-out (FIFO) issuing policy, i.e., older products are issued first to meet demand, which is shown to perform very well in many perishable inventory systems ({e.g., \citet{fries1975optimal2,pierskalla1972optimal}})}. Our contributions in this paper are as follows:

i) {We first present a new approximation algorithm that we call the \textit{marginal-cost dual-balancing policy} for the perishable inventory control problem, and prove that whenever replacing old products with new ones in inventory does not increase the expected total cost, our algorithm has a worst-case performance guarantee of two, i.e., the expected total cost of our policy is at most twice that of an optimal ordering policy.}

ii) {In many perishable inventory systems, the major concern is outdating; and clearly, replacing old products with new ones reduces the chance of products being outdated. {Therefore, the condition that replacing old products with new ones does not increase the expected total cost is very intuitive; however, it is not obvious when this would be ensured to be true {theoretically}.} In that regard, we find that this condition coincides with the optimality of FIFO issuing policy: FIFO being an optimal issuing policy implies that younger products are more preferred to have in inventory than older ones, thus replacing old products with new ones will not increase the expected total cost; and vice versa. Further, given that directly checking the optimality of a FIFO issuing policy may be difficult,  we extend the existing findings in the literature on the optimality of FIFO issuing policy and provide a \textit{necessary and sufficient} condition and several easy-to-check sufficient conditions that ensure FIFO to be an optimal issuing policy.}

{iii) {By ``truncating'' the marginal-cost dual-balancing policy, we present also a more general class of algorithms that we call the \textit{truncated-balancing policy}. In particular, we first prove that a myopic policy under the so-called marginal-cost accounting scheme provides a lower bound on the optimal ordering quantity. Then, we construct the truncated-balancing policy by truncating the marginal-cost dual-balancing policy using the \textit{specific} lower bound we derive and any upper bound on the optimal ordering quantity. We prove that when FIFO is an optimal issuing policy, the truncated-balancing policy also admits a worst-case performance guarantee of two.}

iv) {We further compare our marginal-cost dual-balancing policy with base-stock policies, which are widely studied in the perishable inventory literature. Given that FIFO issuing policy is always optimal under base-stock policies, we show that the expected total cost of the marginal-cost dual-balancing policy is always at most twice that of an optimal base-stock policy.}

v) {Lastly, using real data for the platelet inventory control problem from Hospital Alpha, we conduct extensive computational analyses and show that a) our proposed approximation algorithms perform significantly better than an ``optimal'' policy that does not consider the evolving forecast information, b) the computational performance of our proposed algorithms is substantially better than the theoretical worst-case performance guarantee of two, and c) the truncated-balancing policy has a significant performance improvement over the marginal-cost dual-balancing policy as well as other relevant policies proposed in the literature.}

In the literature, many papers have studied the periodic-review, fixed-lifetime perishable inventory control problem (see reviews by \citet{nahmias1982perishable,nahmias2011perishable} and \citet{karaesmen2011managing}).
The general multi-period lifetime perishable inventory control problem was first studied independently by \citet{fries1975optimal} and \citet{nahmias1975optimal}, who both formulated the problem as a dynamic program (DP) with a state space comprised of inventory levels of different ages. However, the structure of an optimal policy is complicated and finding optimal policies using standard dynamic programming is computationally intractable due to the well-known ``curse of dimensionality''. Therefore, later efforts are mainly focused on heuristic policies. Among the developed heuristic policies, the base-stock policy, under which the total inventory is replenished up to the same level at each period, is particularly popular due to its simplicity and near-optimal numerical performance (e.g., \citet{nahmias1976myopic, cohen1976analysis, chazan1977markovian, nandakumar1993near, cooper2001pathwise, li2009note, chen2014coordinating, zhang2016nonparametric}). Other heuristic policies such as modified base-stock policy (e.g., \citet{broekmeulen2009heuristic}), constant order policy (e.g., \citet{brodheim1975evaluation, deniz2010managing}), and higher-order approximation (\citet{nahmias1977higher}) are also proposed and studied. However, due to the complexity of the perishable inventory control problem, most of these studies assume that demand over time is independently and identically distributed (i.i.d.), and none of the proposed heuristics has a theoretical performance guarantee.

More recently, there is a stream of work focusing on approximation algorithms for stochastic inventory systems under general demand processes. The pioneering work by \citet{levi2007approximation} studies a stochastic inventory control problem for nonperishable products. They show that the proposed dual-balancing policy, which balances the costs of under-ordering and over-ordering under a marginal-cost accounting scheme, has a worst-case performance guarantee of two. This idea has been later extended to many other settings to consider lost sales (\citet{levi20082}), setup costs and capacity constraints (\citet{levi2008approximation, levi2013approximation, shi2014approximation}), remanufacturing (\citet{tao2014approximation}), and perishable products (\citet{chao2016approximation,chao2015approximation}).

{Among these papers that study approximation algorithms in inventory management, \citet{chao2015approximation}, which also considers a perishable inventory control problem with no set-up cost, is the most relevant to ours. In particular, \citet{chao2015approximation} present a proportional-balancing policy and a dual-balancing policy for perishable inventory systems under FIFO issuing policy, and they prove that 1) the proportional-balancing policy has a performance guarantee between two and three for the general case, and 2) the dual-balancing policy has a performance guarantee of two when demand is independent and stochastically non-decreasing over time. While both our study and \citet{chao2015approximation} focus on developing approximation algorithms for perishable inventory systems, our analysis and results are different in the following aspects: i) We present new approximation algorithms for perishable inventory systems (marginal-cost dual balancing policy and truncated-balancing policy) that are different from the ones presented in \citet{chao2015approximation}.} ii) {We tighten the worst-case performance guarantee to exactly two for cases where FIFO is an optimal issuing policy, and show that the condition presented in Chao et al. to ensure a performance guarantee of two (i.e., demand is independent and stochastically non-decreasing over time) is a special case of ours. {Further, we identify several intuitive sufficient conditions that ensure the optimality of FIFO issuing policy, and we present examples where the worst-case performance guarantee of our proposed algorithms is strictly tighter than that presented in \citet{chao2015approximation}} (please see details in $\S$\ref{sec:fifo} and Examples \ref{ex:1} and \ref{ex:2}).} iii) {We further consider truncating the marginal-cost dual-balancing policy using a specific lower bound we derive on the optimal ordering quantity. We remark that this is an important contribution because, unlike the existing results in the nonperishable inventory literature, truncation in the perishable inventory case imposes several new methodological challenges and hence is non-trivial (see also the next paragraph for more details). We show that the truncated-balancing policy also admits a performance guarantee of two, and using real data from a local hospital, we show that it numerically performs much better than the marginal-cost dual-balancing policy and the policies presented in \citet{chao2015approximation}.} iv) {Methodologically, while Chao et al. build their analysis based on algebraic arguments, our analysis is based on two new ideas that we call the \textit{imaginary operation policy} and the \textit{dynamic unit-matching scheme}, respectively (see discussion also in the following paragraph). In particular, the existing worst-case analyses for the nonperishable cases are based on a (static) one-to-one matching between units under two different policies. However, as stated in Chao et al., ``the perishability of products destroys this matching mechanism''. To overcome this challenge, Chao et al. turned to an innovative algebraic approach. On the other hand, the new ideas we propose {(i.e., the imaginary operation policy and the dynamic unit-matching scheme)}, a) significantly simplify the comparison of two different policies and allow us to stay on the track {of} unit matching; b) enable us to reach to a very insightful new result for perishable inventory systems: The worst-case performance guarantee can be tightened to exactly two whenever replacing old products with new ones does not increase the expected total cost (or equivalently when FIFO is an optimal issuing policy); and c) allow us to capture not only perishability but also truncation, which is not considered in Chao et al., but as we show, it significantly improves the computational performance. We believe our ideas are valuable beyond this study and can also be applied to facilitate the analysis for other perishable inventory systems.}}

The idea of truncating a dual-balancing policy with bounds on optimal ordering quantities is first proposed by \citet{hurley2007new}. However, unlike the dual-balancing policy that has been extended to many other settings, the truncated-balancing policy needs a more sophisticated analysis and has only been shown to have a worst-case performance guarantee for the nonperishable backlogging case (\citet{hurley2007new, levi2008approximation}). In this study, we consider a truncated-balancing policy in the perishable inventory setting, where challenges arise from both perishability and the complexity caused by truncation. In particular, first, the analyses for all the nonperishable cases are based on a (static) one-to-one matching between units under two different policies, which rely on the fact that all inventory units will be eventually used to satisfy demand.  However, {this fails to be true for the perishable inventory case, where units can simply outdate without satisfying any demand}. Second, the existing analysis for the truncated-balancing policy in the nonperishable case relies on the base-stock structure of an optimal ordering policy. This also fails to be true for the perishable inventory case, where due to perishability the structure of an optimal policy is complicated and depends on the entire inventory vector. To overcome the above challenges, we introduce two new ideas: 1) a bridging policy that we call the \textit{imaginary operation policy}, under which old products can be replaced with new ones for free so that the inventory vectors under two different policies can be easily compared, and 2) a \textit{dynamic unit-matching scheme}, under which units can be matched and re-matched at different periods so that the shortcomings of the existing static matching approach can be addressed. {As we discussed earlier, these two new ideas together {enable} us to overcome the challenges arising from {both} perishability and truncation, and lead to an insightful {new} result.}

The remainder of this paper is organized as follows. In $\S$\ref{sec:problem}, we present a formal model formulation. In $\S$\ref{sec:accounting}, we present a marginal-cost dual-balancing policy for the perishable inventory control problem.  In $\S$\ref{sec:analysis}, we prove that when FIFO is an optimal issuing policy, our algorithm has a worst-case performance guarantee of two, i.e., the expected total cost of our policy is at most twice that of an optimal ordering policy. We further compare our policy with base-stock policies and show that the expected total cost of our policy is always at most twice that of an optimal base-stock policy. In $\S$\ref{sec:truncated}, we first show that a myopic policy under the marginal-cost accounting scheme provides a lower bound on the optimal ordering quantity; we then present a truncated-balancing policy that also admits a worst-case performance guarantee of two when FIFO is an optimal issuing policy. In $\S$\ref{sec:fifo}, we present a necessary and sufficient condition and several easy-to-check sufficient conditions that ensure the optimality of FIFO issuing policy. Finally, we present computational results based on a platelet inventory control problem in $\S$\ref{sec:numerical}, and draw conclusions in $\S$\ref{sec:conclusion}.

%%%%%%%%%%%%%%%%%%%%%%%%%%%%%%%%%%%%%%%%%%%%%%%%%%%%%
\section{Model Formulation.}
\label{sec:problem}

We study a periodic-review, fixed-lifetime perishable inventory control problem under a general stochastic demand process.

\textbf{Notation}: We consider a product lifetime of $K$ periods and a planning horizon of $T$ periods. Demands over the planning horizon are denoted as $D_1,...,D_T$, which are exogenous random variables with finite means, and can be nonstationary, correlated, and dynamically evolving. As a convention, we generally use capital letters to denote random variables, and lowercase letters to denote their realizations (product lifetime $K$ and planning horizon $T$ are exceptions). At the beginning of each period $t$, there is an information set denoted as $f_t$, which contains the realization of demands $(d_1,...,d_{t-1})$ and possibly some other forecast information available at period $t$, denoted as $(u_1,...,u_t)$. That is, the information set $f_t$ is a specific realization of the random vector $F_t=(D_1,...,D_{t-1},U_1,...,U_t)$. Further, we assume that the conditional joint distribution of future demands $(D_t,...,D_T)$ is known for given $f_t$. Additional notation that describes system states and decision variables is defined as follows:

$X_{k,t}$: the inventory level of age $k$ at the beginning of period $t$, $k=1,...,K-1, t=1,...,T$.

$\textbf{X}_t$: the inventory vector at the beginning of period $t$, i.e., $\textbf{X}_t=(X_{1,t},...,X_{K-1,t})$, $t=1,...,T$. 

$Q_t$: the ordering quantity at period $t, t=1,...,T$. 

$Y_t$: the total inventory level after ordering and before demand realization at period $t$, i.e., $Y_t=\sum\limits_{k=1}^{K-1}X_{k,t}+Q_t$, $t=1,...,T$.

\textbf{System Dynamics}: We define the sequence of events as follows: 1) At the beginning of each period $t=1,...,T$, the $K-1$ dimensional inventory vector $\textbf{X}_t$ and the information set $F_t$ are observed, based on which $Q_t$ products of age 0 are ordered; 2) products ordered arrive instantly with a zero lead time; 3) random demand $D_t$ then occurs during the period, inventory is issued to satisfy demand based on the FIFO rule, and unmet demand is lost (since we assume zero lead time, our results hold equally well for the backlogging case); and 4) at the end of each period, all products in inventory age by 1, and products reaching age $K$ are disposed from the inventory. Let $X_{0,t}=Q_t$. Then, the inventory vector is updated as follows:
$$X_{k,t+1}=\bigg(X_{k-1,t}-\bigg(D_t-\sum\limits_{m=k}^{K-1}X_{m,t}\bigg)^+\bigg)^+, k=1,...,K-1; t=0,...,T-1.$$

Without loss of generality, we assume that the system starts from empty (i.e., zero initial inventory); however, all of our results can be extended to consider arbitrary initial inventory levels.

\textbf{Cost Structure}: At each period, we consider an ordering cost $\hat{c}$ for each unit of product ordered at that period, a shortage penalty $\hat{p}$ for each unit of stock-out, a holding cost $\hat{h}$ for each unit of excess inventory after demand realization, and an outdating cost $\hat{w}$ for each unit of product that is outdated at the end of that period. To eliminate trivial situations, we assume $\hat{p}-\hat{c}\geq0$. We allow negative outdating cost (i.e., positive salvage value) as long as $\hat{w}+\beta\hat{c}\geq0$, where $\beta$ denotes the discount factor. We also consider a salvage value for each unit of product left in inventory at the end of the planning horizon, and for simplicity we assume it is equal to the ordering cost $\hat{c}$ (our results can be easily extended to consider any salvage value $\hat{v}$ as long as $\hat{w}+\beta\hat{v}\geq0$). 

\textbf{Optimality Criterion}: At each period $t$, given the inventory vector $\textbf{x}_t$ and the information set $f_t$, an ordering \textit{decision rule} is a function from the set of all possible $(\textbf{x}_t,f_t)$ to the set of all possible $q_t$; and an ordering \textit{policy} is a collection of ordering decision rules at all periods. Let $\pi$ denote any given ordering policy. Then, the total cost under policy $\pi$ over the planning horizon is:
\begin{equation*}
\hat{\mathscr{C}}(\pi)= \sum\limits_{t=1}^T\beta^{t-1}\big(\hat{c}Q_t^{\pi}+\hat{p}(D_t-Y_t^{\pi})^++\hat{h}(Y_t^{\pi}-D_t)^++\hat{w}(X_{K-1,t}^{\pi}-D_t)^+\big)-\beta^{T}\hat{c}\sum\limits_{k=1}^{K-1}X_{k,T+1}^{\pi},
\end{equation*}
where $\textbf{X}_t^{\pi}$ and $Q_t^{\pi}$ denote the inventory vector and the ordering quantity at period $t$ under policy $\pi$, respectively. Then, our problem is to find an optimal ordering policy $OPT$ such that $OPT\in \arg\min\limits_{\pi}\mathrm{E}[\hat{\mathscr{C}}(\pi)]$.

%%%%%%%%%%%%%%%%%%%%%%%%%%%%%%%%%%%%%%%%%%%%%%%%%%%%%%
\section{Marginal-Cost Dual-Balancing Policy.}
\label{sec:accounting}

In this section, we first introduce a cost transformation to eliminate ordering cost in $\S$\ref{sec:cost}. We then present a marginal-cost accounting scheme for the perishable inventory setting in $\S$\ref{sec:marginal}. Finally, we present our algorithm in $\S$\ref{sec:algorithm}. Unless presented in the main text, the proofs of all analytical results are included in the Appendix.

\subsection{Cost Transformation.}
\label{sec:cost}

To apply the marginal-cost accounting scheme which we present in $\S$\ref{sec:marginal}, we first need to construct an equivalent problem with a zero ordering cost.

Define the cost parameters for the transformed problem as: $c=0, p=\hat{p}-\hat{c}, h=\hat{h}+(1-\beta)\hat{c}$, and $w=\hat{w}+\beta \hat{c}$. {Since we assume $\hat{p}-\hat{c}\geq0$ and $\hat{w}+\beta\hat{c}\geq0$, all the transformed costs are nonnegative}. Then, for a given policy $\pi$, the total cost of the transformed problem is:
\begin{equation*}
\mathscr{C}(\pi)= \sum\limits_{t=1}^T\beta^{t-1}\big(p(D_t-Y_t^{\pi})^++h(Y_t^{\pi}-D_t)^++w(X_{K-1,t}^{\pi}-D_t)^+\big).
\end{equation*}

In the following lemma, we show that the difference between the total costs of the original and transformed problems is independent of policy $\pi$, which implies that the two problems are equivalent in the sense that they have the same set of optimal ordering policies.

\begin{lemma}
\label{lem:c1c2}
For any policy $\pi$, $\hat{\mathscr{C}}(\pi)-\mathscr{C}(\pi)=\sum\limits_{t=1}^T\beta^{t-1}\hat{c}D_t$, with probability one.
\end{lemma}

\subsection{Marginal-Cost Accounting Scheme.}
\label{sec:marginal}

Unlike traditional methods which assign each period all costs that occur at this period, the marginal-cost accounting scheme, introduced by \citet{levi2007approximation}, assigns each period all costs that are caused by the decision made at this period. For example, a unit ordered at period $t$ may stay in the system for multiple periods, thus holding costs may be charged for this unit for multiple periods; under the marginal-cost accounting scheme, all these holding costs are assigned to period $t$. We now present the marginal-cost accounting scheme for the perishable inventory setting.

\textbf{Marginal Shortage Penalty}: Since inventory can be replenished with a zero lead time, the marginal shortage penalty at each period is simply defined as the shortage penalty that occurs at this period. For $t=1,...,T$, given $\textbf{x}_t, f_t$ and $q_t$, let $P_t(\textbf{x}_t,f_t,q_t)$ denote the expected marginal shortage penalty at period $t$. Then, we have:
$$P_t(\textbf{x}_t,f_t,q_t)\mathrel{\mathop:}=\beta^{t-1}p\mathrm{E}[(D_t-y_t)^+|f_t].$$

\textbf{Marginal Holding Cost}: For $t=1,...,T$, given $\textbf{x}_t$, $f_t$ and $q_t$, let $H_t(\textbf{x}_t,f_t,q_t)$ denote the expected marginal holding cost at period $t$, which is defined as the sum of all expected holding costs charged for units ordered at period $t$. In the perishable inventory setting, since units in inventory may become outdated without satisfying any demand, the future holding costs charged for $q_t$ depend on the entire inventory vector $\textbf{x}_t$. Thus, similar to \citet{nahmias1975optimal}, we let $A_{0,t}=0$, and for $k=1,...,K-1$, let $A_{k,t}$ be the total demand over periods $t, ...,t+k-1$ that cannot be satisfied by the inventory of $(x_{K-k,t},...,x_{K-1,t})$, i.e., the inventory that would have been outdated by the end of period $t+k-1$. Then:
$$A_{k,t}=(A_{k-1,t}+D_{t+k-1}-x_{K-k,t})^+, k=1,...,K-1.$$
Thus, for $k=0,...,K-1$, $(A_{k,t}+D_{t+k}-\sum\limits_{m=1}^{K-k-1}x_{m,t})^+$ represents the total demand over periods $t, ...,t+k$ that cannot be satisfied by the inventory of $\textbf{x}_t$, and $(q_t-(A_{k,t}+D_{t+k}-\sum\limits_{m=1}^{K-k-1}x_{m,t})^+)^+$ represents the amount of $q_t$ left in inventory at the end of period $t+k$. Then, we have:
$$H_t(\textbf{x}_t,f_t,q_t)\mathrel{\mathop:} =\sum\limits_{k=0}^{K-1}\beta^{t+k-1}h\mathrm{E}\bigg[(q_t-(A_{k,t}+D_{t+k}-\sum\limits_{m=1}^{K-k-1}x_{m,t})^+)^+\bigg|f_t\bigg],$$
where the sum over $k$ is defined up to $T-t$ when $t+K-1\geq T$.

\textbf{Marginal Outdating Cost}: For $t=1,...,T$, given $\textbf{x}_t$, $f_t$ and $q_t$, let $W_t(\textbf{x}_t,f_t,q_t)$ denote the expected marginal outdating cost at period $t$, which is defined as the sum of all expected outdating costs charged for units ordered at period $t$, i.e., the expected outdating costs that occur at period $t+K-1$. Note that units ordered at periods $T-K+2,...,T$ will not be outdated within the planning horizon, thus we simply define $W_t(\textbf{x}_t,f_t,q_t)=0$ for $t=T-K+2,...,T$. For $t\leq T-K+1$, $(q_t-A_{K-1,t}-D_{t+K-1})^+$ represents the amount of $q_t$ that will be outdated at the end of period $t+K-1$. Then, we have:
$$W_t(\textbf{x}_t,f_t,q_t)\mathrel{\mathop:}=\beta^{t+K-1}w\mathrm{E}[(q_t-A_{K-1,t}-D_{t+K-1})^+|f_t].$$

For a given policy $\pi$, let $P_t^{\pi}, H_t^{\pi}$ and $W_t^{\pi}$ denote the corresponding marginal shortage penalty, holding and outdating costs at period $t$, respectively. Under a given policy $\pi$, $\textbf{x}^{\pi}_t$ and $q_t^{\pi}$ are both known for given $f_t$. Then, $\mathrm{E}[P_t^{\pi}|f_t]=P_t(\textbf{x}_t^{\pi},f_t,q_t^{\pi})$, $\mathrm{E}[H_t^{\pi}|f_t]=H_t(\textbf{x}_t^{\pi},f_t,q_t^{\pi})$, and $\mathrm{E}[W_t^{\pi}|f_t]=W_t(\textbf{x}_t^{\pi},f_t,q_t^{\pi})$. Since the system starts from zero inventory, we have $\mathscr{C}(\pi)=\sum\limits_{t=1}^T(P_t^{\pi}+H_t^{\pi}+W_t^{\pi}).$ 

\subsection{Algorithm.}
\label{sec:algorithm}
Now we present our first algorithm based on the marginal-cost accounting scheme presented above. Clearly, the expected marginal shortage penalty $P_t(\textbf{x}_t,f_t,q_t)$ occurs due to under-ordering, while the expected marginal holding and outdating costs $H_t(\textbf{x}_t,f_t,q_t)$ and $W_t(\textbf{x}_t,f_t,q_t)$ occur due to over-ordering. Therefore, we define the \textit{marginal-cost dual-balancing policy} (denoted as $B$) as to balance the expected marginal shortage penalty against the sum of the expected marginal holding and outdating costs. More specifically, at each period $t$, given $\textbf{x}_t$ and $f_t$, the marginal-cost dual-balancing ordering quantity $q_t^B$ (for simplicity, we also call it the balancing ordering quantity in the following text) is defined as the solution to the following equation:
\begin{equation}
\label{eqn:balancing}
P_t(\textbf{x}_t,f_t,q_t)=H_t(\textbf{x}_t,f_t,q_t)+W_t(\textbf{x}_t,f_t,q_t).
\end{equation}

Note that the existence of the balancing ordering quantity $q_t^B$ is guaranteed, because at any period $t$, given $\textbf{x}_t$ and $f_t$, $P_t(\textbf{x}_t,f_t,q_t)$ is non-increasing in $q_t$; when $q_t=0$, $P_t(\textbf{x}_t,f_t,q_t)$ is nonnegative, and when $q_t$ goes to infinity, $P_t(\textbf{x}_t,f_t,q_t)$ goes to zero (since demand has a finite mean). In contrast, $H_t(\textbf{x}_t,f_t,q_t)$ and $W_t(\textbf{x}_t,f_t,q_t)$ are non-decreasing in $q_t$; when $q_t=0$, $H_t(\textbf{x}_t,f_t,q_t)=W_t(\textbf{x}_t,f_t,q_t)=0$, and when $q_t$ goes to infinity, both $H_t(\textbf{x}_t,f_t,q_t)$ and $W_t(\textbf{x}_t,f_t,q_t)$ go to infinity. Therefore, $q_t^B$ is guaranteed to exist when we allow fractional ordering quantities. The algorithm can be easily extended to consider discrete ordering quantities following a similar argument as in \citet{levi2007approximation}.

We also remark that our marginal-cost dual-balancing policy is different from the dual-balancing policy defined in \citet{chao2015approximation}. In particular, while our marginal-cost dual-balancing policy balances the marginal shortage penalty against the sum of the marginal holding and outdating costs, the dual-balancing policy in \citet{chao2015approximation} balances the marginal shortage penalty against the marginal outdating cost plus the holding cost that occurs at period $t$, i.e., the marginal holding cost $H_t(\textbf{x}_t,f_t,q_t)$ in Equation (\ref{eqn:balancing}) is replaced by $\beta^{t-1}h\mathrm{E}[(y_t-D_t)^+|f_t]$ in \citet{chao2015approximation}.

%%%%%%%%%%%%%%%%%%%%%%%%%%%%%%%%%%%%%%%%%%%%%%%%%%%%%%
\section{Worst-Case Analysis.}
\label{sec:analysis}

In this section, we first build a bridging policy in $\S$\ref{sec:im}. Then, in $\S$\ref{sec:matching}, we construct a new unit-matching scheme that (dynamically) matches units under two different policies on a one-to-one correspondence. Based on these results, we show in $\S$\ref{sec:worst} that when FIFO is an optimal issuing policy, our marginal-cost dual-balancing policy has a worst-case performance guarantee of two, i.e., the expected total cost of our policy is at most twice that of an optimal ordering policy. Finally, in $\S$\ref{sec:base}, we compare our policy with an optimal base-stock policy, and show that the expected total cost of our policy is always at most twice that of an optimal base-stock policy.

%%%%%%%%%%%%%%%%%%%%%%%%%%%%%%%%%%%%%%%%%%%%%%%%%%%%%
\subsection{A Bridging Policy: Imaginary Operation Policy.}
\label{sec:im}

By Lemma \ref{lem:c1c2}, we know that $\hat{\mathscr{C}}(\pi)-\mathscr{C}(\pi)$ is nonnegative and independent of policy $\pi$. Therefore, to show that the expected total cost of the marginal-cost dual-balancing policy is at most twice that of an optimal ordering policy (i.e., $\mathrm{E}[\hat{\mathscr{C}}(B)]\leq 2\mathrm{E}[\hat{\mathscr{C}}(OPT)]$), it is sufficient to show that $\mathrm{E}[\mathscr{C}(B)]\leq 2\mathrm{E}[\mathscr{C}(OPT)].$

However, due to the partially ordered nature of multi-dimensional inventory vectors, it is difficult to directly compare the costs under policies $B$ and $OPT$ . Therefore, we next propose a bridging policy that we call the \textit{imaginary operation policy} (denoted as $IM$), which allows us to properly modify the inventory vectors so that the inventory vectors under two different policies become completely ordered, and the respective costs can be easily compared. We then show $\mathrm{E}[\mathscr{C}(IM)]\leq \mathrm{E}[\mathscr{C}(OPT)]$ (Lemma \ref{lem:optim}) and $\mathrm{E}[\mathscr{C}(B)]\leq 2\mathrm{E}[\mathscr{C}(IM)]$ (Lemma \ref{lem:imb}), respectively, which leads to our main result $\mathrm{E}[\mathscr{C}(B)]\leq 2\mathrm{E}[\mathscr{C}(OPT)]$ (Theorem \ref{optb}).

Policy $IM$ is constructed as follows: At each period $t$, given $\textbf{x}_t$ and $f_t$, let the system under policy $IM$ follow an optimal ordering policy.\footnote{Throughout the paper, we refer following an optimal ordering policy to implementing an optimal decision rule at each period given the system state, instead of copying the ordering quantity from the system under policy $OPT$.} What differentiates policies $IM$ and $OPT$ is that under policy $IM$, {at each period after ordering and before demand realization}, products in the inventory vector can be ``moved'' from older positions to the position of age 0, i.e., old products can be replaced with new ones for free. Note that since the inventory vectors under policies $IM$ and $OPT$ may be different at each period, the actual ordering quantities under the two policies can also be different. 

At each period $t$, let $y_t^{B}$ and $y_t^{IM}$ be the total inventory levels after ordering under policies $B$ and $IM$, respectively (note that once the rules of movements under policy $IM$ for periods $1,...,t-1$ are defined, $y_t^{IM}$ is well-defined). Then, we partition the set of decision epochs $\{1,...,T\}$ into the following two subsets:
$$\mathscr{T}_P=\{t: y_t^B\geq y_t^{IM}\}, \mathscr{T}_H=\{t: y_t^B< y_t^{IM}\}.$$

The main objective of constructing policy $IM$ is to bound the total shortage penalty of policy $B$ at each period $t\in \mathscr{T}_P$ and the total holding and outdating costs of policy $B$ charged for the units ordered at each period $t\in \mathscr{T}_H$. {Since we have $y_t^B\geq y_t^{IM}, \forall t\in \mathscr{T}_P$, the total shortage penalty of policy $B$ at $t\in \mathscr{T}_P$ can be easily bounded. Therefore, unit movements are only needed at $t\in \mathscr{T}_H$.} The rules of movements are defined as follows, and an illustrative example is provided at the end of this subsection.

{Let $\mathscr{T}_H=\{\tau_1,...,\tau_n\}$, where $\tau_1<...<\tau_n$. At the beginning of period $\tau_1$ after ordering (and before demand realization), we simply move all units in the inventory vector under policy $IM$ to age 0.}

{At $\tau_2$, units ordered at $\tau_1$ become of age $\tau_2-\tau_1$. Since all units are moved to age 0 at $\tau_1$, the inventory under policy $IM$ is consumed (used to satisfy demand or outdated) no faster than that under policy $B$. Also, we have $y_{\tau_1}^{IM}> y_{\tau_1}^{B}$ at $\tau_1$. Then at $\tau_2$, the total inventory of age greater than or equal to age $\tau_2-\tau_1$ under policy $IM$ is no less than that under policy $B$, i.e., $\sum\limits_{k=\tau_2-\tau_1}^{K-1}x_{k,\tau_2}^{IM}\geq \sum\limits_{k=\tau_2-\tau_1}^{K-1}x_{k,\tau_2}^{B}$. Therefore, at the beginning of period $\tau_2$ after ordering, we first move all units of age strictly less than $\tau_2-\tau_1$ under policy $IM$ to age 0 such that there are only positive inventory of age 0 and $\tau_2-\tau_1$ under policy $IM$. We then move some units of age equal to $\tau_2-\tau_1$ under policy $IM$ to age 0 such that $\sum\limits_{k=\tau_2-\tau_1}^{K-1}x_{k,\tau_2}^{IM}= \sum\limits_{k=\tau_2-\tau_1}^{K-1}x_{k,\tau_2}^{B}$.}

{Similarly, for any $i\geq 2$, at the beginning of period $\tau_i$ after ordering, we first move all units of age strictly less than $\tau_i-\tau_{i-1}$ under policy $IM$ to age 0, and then for each $j=1,...,i-1$, move some units of age equal to $\tau_i-\tau_j$ under policy $IM$ to age 0 such that after all the movements, we have:}

{(i) There are only positive inventory of age $0, \tau_i-\tau_{i-1},...,\tau_i-\tau_1$ under policy $IM$.}

{(ii) For $j=1,...,i-1$, the total inventory of age greater than or equal to age $\tau_i-\tau_j$ under policies $IM$ and $B$ are the same, i.e.,}

\begin{equation}
\label{eqn:move1}
{\sum\limits_{k=\tau_i-\tau_j}^{K-1}x_{k,\tau_i}^{IM}=\sum\limits_{k=\tau_i-\tau_j}^{K-1}x_{k,\tau_i}^{B}, j=1,...,i-1.}
\end{equation}

{Based on the rules of movements defined above, we are ensured to have that: $\forall \tau_i\in\mathscr{T}_H$, after the movements of units at $\tau_i$, the inventory vector under policy $IM$ is ``younger'' than that under policy $B$, i.e., for $k=1,...,K-1$, policy $IM$ has no more inventory of age greater than or equal to $k$.}

\begin{lemma}
\label{lem:im}
$\forall \tau_i\in\mathscr{T}_H$, after the movements of units at $\tau_i$, we have:
\begin{equation}
\sum\limits_{m=k}^{K-1}x_{m,\tau_i}^{IM}\leq \sum\limits_{m=k}^{K-1}x_{m,\tau_i}^{B}, k=1,...,K-1.
\label{Eq. ine}
\end{equation}
\end{lemma}

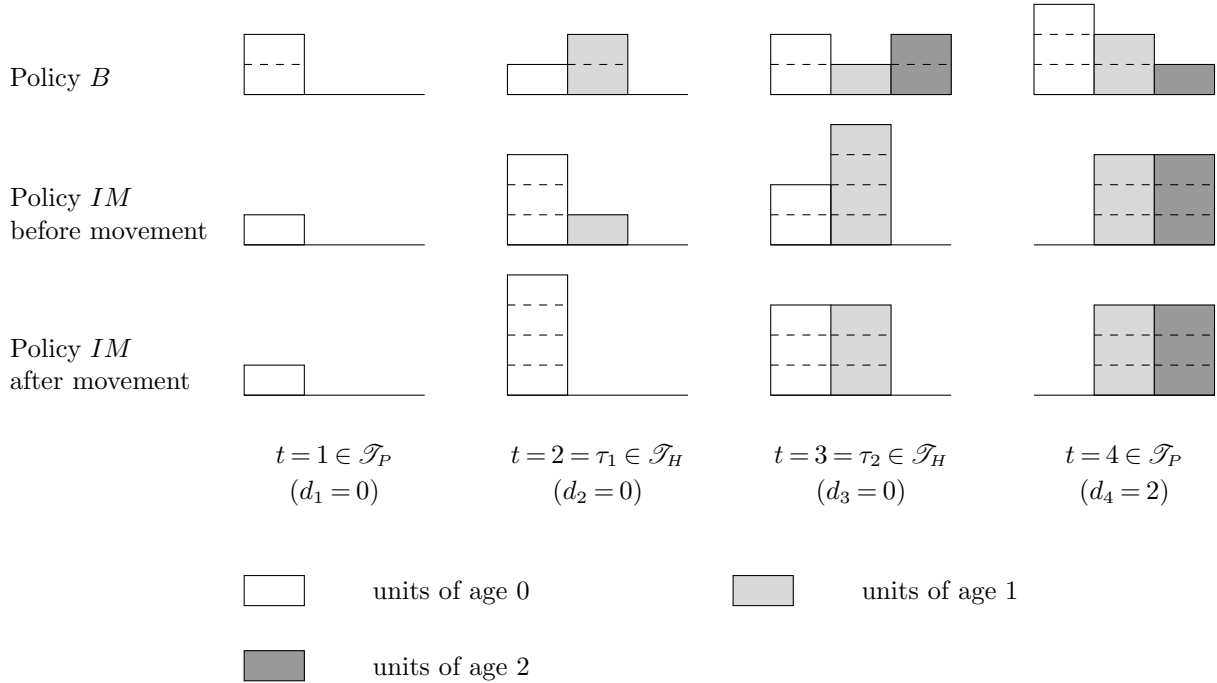
\begin{figure}[h]
\begin{tikzpicture}
\small
\draw (2,0) rectangle (2.8,0.4);
\draw (2,0) rectangle (4.4,0);
\draw (2,2) rectangle (2.8,2.4);
\draw (2,2) rectangle (4.4,2);
\draw (2,4) rectangle (2.8,4.8);
\draw (2,4) rectangle (4.4,4);
\draw[dashed] (2,4.4) -- (2.8,4.4);
\draw (-0.445,4.2) node {Policy $B$};
\draw (-0.3,2.6) node {Policy $IM$};
\draw (0.2,2.2) node {before movement};
\draw (-0.3,0.6) node {Policy $IM$};
\draw (0.08,0.2) node {after movement};
\draw (3.2,-0.8) node {$t=1\in\mathscr{T}_P$};
\draw (3.2,-1.3) node {$(d_1=0)$};

\draw (5.5,0) rectangle (6.3,1.6);
\draw[dashed] (5.5,0.4) -- (6.3,0.4);
\draw[dashed] (5.5,0.8) -- (6.3,0.8);
\draw[dashed] (5.5,1.2) -- (6.3,1.2);
\draw (5.5,0) rectangle (7.9,0);
\draw (5.5,2) rectangle (6.3,3.2);
\draw[dashed] (5.5,2.4) -- (6.3,2.4);
\draw[dashed] (5.5,2.8) -- (6.3,2.8);
\fill[black!15] (6.3,2) -- (6.3,2.4) -- (7.1,2.4) -- (7.1,2) -- cycle;
\draw (6.3,2) rectangle (7.1,2.4);
\draw (5.5,2) rectangle (7.9,2);
\draw (5.5,4) rectangle (6.3,4.4);
\fill[black!15] (6.3,4) -- (6.3,4.8) -- (7.1,4.8) -- (7.1,4) -- cycle;
\draw (6.3,4) rectangle (7.1,4.8);
\draw[dashed] (6.3,4.4) -- (7.1,4.4);
\draw (5.5,4) rectangle (7.9,4);
\draw (6.7,-0.8) node {$t=2=\tau_1\in\mathscr{T}_H$};
\draw (6.7,-1.3) node {$(d_2=0)$};

\draw (9,0) rectangle (9.8,1.2);
\draw[dashed] (9,0.4) -- (9.8,0.4);
\draw[dashed] (9,0.8) -- (9.8,0.8);
\fill[black!15] (9.8,0) -- (9.8,1.2) -- (10.6,1.2) -- (10.6,0) -- cycle;
\draw (9.8,0) rectangle (10.6,1.2);
\draw[dashed] (9.8,0.4) -- (10.6,0.4);
\draw[dashed] (9.8,0.8) -- (10.6,0.8);
\draw (9,0) rectangle (11.4,0);
\fill[black!15] (9.8,2) -- (9.8,3.6) -- (10.6,3.6) -- (10.6,2) -- cycle;
\draw (9.8,2) rectangle (10.6,3.6);
\draw (9,2) rectangle (9.8,2.8);
\draw[dashed] (9,2.4) -- (9.8,2.4);
\draw[dashed] (9.8,2.4) -- (10.6,2.4);
\draw[dashed] (9.8,2.8) -- (10.6,2.8);
\draw[dashed] (9.8,3.2) -- (10.6,3.2);
\draw (9,2) rectangle (11.4,2);
\draw (9,4) rectangle (9.8,4.8);
\draw[dashed] (9,4.4) -- (9.8,4.4);
\fill[black!15] (9.8,4) -- (9.8,4.4) -- (10.6,4.4) -- (10.6,4) -- cycle;
\draw (9.8,4) rectangle (10.6,4.4);
\fill[black!40] (10.6,4) -- (10.6,4.8) -- (11.4,4.8) -- (11.4,4) -- cycle;
\draw (10.6,4) rectangle (11.4,4.8);
\draw[dashed] (10.6,4.4) -- (11.4,4.4);
\draw (9,4) rectangle (11.4,4);
\draw (10.2,-0.8) node {$t=3=\tau_2\in\mathscr{T}_{H}$};
\draw (10.2,-1.3) node {$(d_3=0)$};

\fill[black!15] (13.3,0) -- (13.3,1.2) -- (14.1,1.2) -- (14.1,0) -- cycle;
\draw (13.3,0) rectangle (14.1,1.2);
\draw[dashed] (13.3,0.4) -- (14.1,0.4);
\draw[dashed] (13.3,0.8) -- (14.1,0.8);
\fill[black!40] (14.1,0) -- (14.1,1.2) -- (14.9,1.2) -- (14.9,0) -- cycle;
\draw (14.1,0) rectangle (14.9,1.2);
\draw[dashed] (14.1,0.4) -- (14.9,0.4);
\draw[dashed] (14.1,0.8) -- (14.9,0.8);
\draw (12.5,0) rectangle (14.9,0);
\fill[black!15] (13.3,2) -- (13.3,3.2) -- (14.1,3.2) -- (14.1,2) -- cycle;
\draw (13.3,2) rectangle (14.1,3.2);
\draw[dashed] (13.3,2.4) -- (14.1,2.4);
\draw[dashed] (13.3,2.8) -- (14.1,2.8);
\fill[black!40] (14.1,2) -- (14.1,3.2) -- (14.9,3.2) -- (14.9,2) -- cycle;
\draw (14.1,2) rectangle (14.9,3.2);
\draw[dashed] (14.1,2.4) -- (14.9,2.4);
\draw[dashed] (14.1,2.8) -- (14.9,2.8);
\draw (12.5,2) rectangle (14.9,2);
\draw (12.5,4) rectangle (13.3,5.2);
\draw[dashed] (12.5,4.4) -- (13.3,4.4);
\draw[dashed] (12.5,4.8) -- (13.3,4.8);
\fill[black!15] (13.3,4) -- (13.3,4.8) -- (14.1,4.8) -- (14.1,4) -- cycle;
\draw (13.3,4) rectangle (14.1,4.8);
\draw[dashed] (13.3,4.4) -- (14.1,4.4);
\fill[black!40] (14.1,4) -- (14.1,4.4) -- (14.9,4.4) -- (14.9,4) -- cycle;
\draw (14.1,4) rectangle (14.9,4.4);
\draw (12.5,4) rectangle (14.9,4);
\draw (13.7,-0.8) node {$t=4\in\mathscr{T}_P$};
\draw (13.7,-1.3) node {$(d_4=2)$};

\draw (2,-2.8) rectangle (2.8,-2.4);
\draw (4.75,-2.63) node {units of age 0};
\fill[black!15] (8.5,-2.8) -- (9.3,-2.8) -- (9.3,-2.4) -- (8.5,-2.4) -- cycle;
\draw (8.5,-2.8) rectangle (9.3,-2.4);
\draw (11.25,-2.63) node {units of age 1};
\fill[black!40] (2,-3.8) -- (2.8,-3.8) -- (2.8,-3.4) -- (2,-3.4) -- cycle;
\draw (2,-3.8) rectangle (2.8,-3.4);
\draw (4.75,-3.63) node {units of age 2};
\end{tikzpicture}
\caption{An illustrative example to show the imaginary operation policy ($IM$)}
\label{figure:Figure1}
\end{figure}

An illustrative example describing the rules of movements is presented in Figure \ref{figure:Figure1}. In this example, we have product lifetime of $K=3$ periods, and planning horizon of $T=4$ periods. Consider a given sample path where $d_1=d_2=d_3=0$ and $d_4=2$. At the beginning of period $t=1$, assume $q_1^B=2$ and $q_1^{IM}=1$. Then, $y_1^{B}=2>1=y_1^{IM}$, thus $t=1\in \mathscr{T}_P$ and no movements are performed at this period. At the beginning of period $t=2$, assume $q_2^B=1$ and $q_2^{IM}=3$. Then, $y_2^{B}=3<4=y_2^{IM}$, thus $t=2=\tau_1\in \mathscr{T}_H$, and we move all units under policy $IM$ to age 0. At the beginning of period $t=3$, assume $q_3^{B}=q_3^{IM}=2$. Then, $y_3^{B}=5<6=y_3^{IM}$, thus $t=3=\tau_2\in \mathscr{T}_{H}$. The unit ordered at $\tau_1$ under policy $B$ is now of age $\tau_2-\tau_1=1$, therefore we move one unit of age 1 under policy $IM$ to age 0 such that the amount of units of age greater than or equal to 1 under policies $B$ and $IM$ are equal. At the beginning of period $t=4$, assume $q_4^B=3$ and $q_4^{IM}=0$. Then, $y_4^{B}=6=y_4^{IM}$, thus $t=4\in \mathscr{T}_P$ and no movements are performed at this period.

%%%%%%%%%%%%%%%%%%%%%%%%%%%%%%%%%%%%%%%%%%%%%%%%%%%%%
\subsection{A Dynamic Unit-Matching Scheme.}
\label{sec:matching}

Based on the imaginary operation policy $(IM)$ we constructed above, we now introduce a new unit-matching scheme that matches inventory units under policies $B$ and $IM$, which plays a key role in the comparison of $\mathscr{C}(B)$ and $\mathscr{C}(IM)$. In particular, our objective is to match the units ordered at each period $t\in\mathscr{T}_H$ under policy $B$ to units under policy $IM$ on a one-to-one correspondence, such that a matched unit under policy $B$ stays in inventory no longer than the corresponding unit under policy $IM$. This way, the total holding and outdating costs charged for the units ordered at $t\in\mathscr{T}_H$ under policy $B$ can be bounded by the total holding and outdating costs under policy $IM$. 

The idea of examining inventory and demand at a unit level is first proposed by \citet{muharremoglu2008single}, and is first applied to prove worst-case performance guarantee by \citet{levi2007approximation}, where units under two policies are matched on a one-to-one correspondence. Similar arguments are also used in all the subsequent studies on approximation algorithms for nonperishable inventory systems (\citet{levi20082, levi2008approximation, levi2013approximation, shi2014approximation, tao2014approximation}). However, in these studies, the matching of inventory units is static in the sense that once a pair of units under two policies are matched at some period, the matching is permanent. This approach relies on the assumption that all units ordered will be eventually used to satisfy demand, and a pair of units, once matched, will be used to satisfy the same unit of demand. However, this fails to be true in the perishable inventory setting, where units in inventory may simply outdate without satisfying any demand. To address this complication, we introduce a new matching scheme that we call the \textit{dynamic unit-matching scheme}, under which, a unit ordered at $t\in\mathscr{T}_H$ under policy $B$ can be matched and then re-matched to a new unit under policy $IM$. The rules of matchings are defined as follows, and an illustrative example is provided at the end of this subsection.

{Recall that $\mathscr{T}_H=\{\tau_1,...,\tau_n\}$, where $\tau_1<...<\tau_n$}. At the beginning of period $\tau_1$, after the movements of units under policy $IM$ (based on the rules described in $\S$\ref{sec:im}), we assign indices from 1 to $y_{\tau_1}^{B}$ for units under policy $B$, and assign indices from 1 to $y_{\tau_1}^{IM}$ for units under policy $IM$,\footnote{For continuous demands and ordering quantities, indices are defined continuously from 0 to $y_{\tau_1}^{B}$ and $y_{\tau_1}^{IM}$.} where $y_{\tau_1}^B< y_{\tau_1}^{IM}$. Older units are assigned smaller indices, and units of the same age are sorted in an arbitrary sequence and assigned indices accordingly. Then, we \textit{temporarily} match each unit ordered at $\tau_1$ under policy $B$ to the unit with the same index under policy $IM$. Clearly, each pair of temporarily matched units have the same age (of age 0).

Then, at the beginning of period $\tau_2$, consider the following three cases. First, if a temporarily matched unit under policy $B$ has been used to satisfy demand, there must exist a unit under policy $IM$ that is also used to satisfy the same unit of demand. This is because $y_{\tau_1}^B< y_{\tau_1}^{IM}$ and the inventory under policy $IM$ is consumed no faster than that under policy $B$ (due to Inequality (\ref{Eq. ine})). Therefore, for a temporarily matched unit under policy $B$ that has been used to satisfy demand, we re-match it to the unit under policy $IM$ that is used to satisfy the same unit of demand, and we set this matching to be \textit{permanent}. 

Second, if a temporarily matched unit under policy $B$ has been outdated, its last temporarily matched unit under policy $IM$ must also have been outdated (since they have the same age). We set this matching also to be permanent. 

Third, for units that are still in inventory at the beginning of period $\tau_2$, we re-define the indices and re-match them to new units under policy $IM$. In particular, we assign indices from 1 to $y_{\tau_2}^{B}$ to units under policy $B$ and assign indices from 1 to $y_{\tau_2}^{IM}$ to units under policy $IM$. Then, we re-match (still temporarily) all previously temporarily matched units (now of age $\tau_2-\tau_1$) and all units ordered at $\tau_2$ (now of age 0) under policy $B$ to units with the same indices under policy $IM$. Since after the movements of units at $\tau_2$, there are only positive inventory of age 0 and $\tau_2-\tau_1$ under policy $IM$ and Equation (\ref{eqn:move1}) holds, each pair of temporarily matched units must have the same age (either 0 or $\tau_2-\tau_1$).

Continuing in this manner, all units ordered at $t\in\mathscr{T}_H$ under policy $B$ are ultimately permanently matched to certain units under policy $IM$ if they are used to satisfy demand or outdated. For units that are still in inventory at the end of the planning horizon, their last temporarily matched units under policy $IM$ must also be in inventory. This is because after the movements at the last period in $\mathscr{T}_H$, the inventory vector under policy $IM$ is ``younger'' than that under policy $B$ (i.e., Inequality (\ref{Eq. ine})), which ensures that a unit under policy $IM$ is consumed no earlier than the unit with the same index under policy $B$. We then set these matchings also to be permanent. Note that by construction, there are no overlaps in the permanent matchings. Further, as we will show in Lemma \ref{lem:imb}, any matched unit under policy $B$ stays in inventory no longer than its permanently matched unit under policy $IM$.

\begin{figure}[h]
\begin{tikzpicture}
\small
\draw (2,0) rectangle (2.8,0.4);
\draw (2,0) rectangle (4.4,0);
\draw (2,2) rectangle (2.8,2.4);
\draw (2,2) rectangle (4.4,2);
\draw (2,4) rectangle (2.8,4.8);
\draw (2,4) rectangle (4.4,4);
\draw[dashed] (2,4.4) -- (2.8,4.4);
\draw (-0.445,4.2) node {Policy $B$};
\draw (-0.3,2.6) node {Policy $IM$};
\draw (0.2,2.2) node {before movement};
\draw (-0.3,0.6) node {Policy $IM$};
\draw (0.08,0.2) node {after movement};
\draw (3.2,-0.8) node {$t=1\in\mathscr{T}_P$};
\draw (3.2,-1.3) node {$(d_1=0)$};

\draw (5.5,0) rectangle (6.3,1.6);
\draw[dashed] (5.5,0.4) -- (6.3,0.4);
\draw[dashed] (5.5,0.8) -- (6.3,0.8);
\draw[dashed] (5.5,1.2) -- (6.3,1.2);
\draw (5.5,0) rectangle (7.9,0);
\draw (5.5,2) rectangle (6.3,3.2);
\draw[dashed] (5.5,2.4) -- (6.3,2.4);
\draw[dashed] (5.5,2.8) -- (6.3,2.8);
\fill[black!15] (6.3,2) -- (6.3,2.4) -- (7.1,2.4) -- (7.1,2) -- cycle;
\draw (6.3,2) rectangle (7.1,2.4);
\draw (5.5,2) rectangle (7.9,2);
\draw[pattern=north west lines, pattern color=black!60] (5.5,4) rectangle (6.3,4.4);
\draw (5.5,4) rectangle (6.3,4.4);
\fill[black!15] (6.3,4) -- (6.3,4.8) -- (7.1,4.8) -- (7.1,4) -- cycle;
\draw (6.3,4) rectangle (7.1,4.8);
\draw[dashed] (6.3,4.4) -- (7.1,4.4);
\draw (5.5,4) rectangle (7.9,4);
\draw (6.7,-0.8) node {$t=2=\tau_1\in\mathscr{T}_H$};
\draw (6.7,-1.3) node {$(d_2=0)$};
\draw (6.7,4.6) node {1};
\draw (6.7,4.2) node {2};
\draw (5.9,4.2) node {3};
\draw (5.9,1.4) node {1};
\draw (5.9,1) node {2};
\draw (5.9,0.6) node {3};
\draw (5.9,0.2) node {4};

\draw (9,0) rectangle (9.8,1.2);
\draw[dashed] (9,0.4) -- (9.8,0.4);
\draw[dashed] (9,0.8) -- (9.8,0.8);
\fill[black!15] (9.8,0) -- (9.8,1.2) -- (10.6,1.2) -- (10.6,0) -- cycle;
\draw (9.8,0) rectangle (10.6,1.2);
\draw[dashed] (9.8,0.4) -- (10.6,0.4);
\draw[dashed] (9.8,0.8) -- (10.6,0.8);
\draw (9,0) rectangle (11.4,0);
\fill[black!15] (9.8,2) -- (9.8,3.6) -- (10.6,3.6) -- (10.6,2) -- cycle;
\draw (9,2) rectangle (9.8,2.8);
\draw[dashed] (9,2.4) -- (9.8,2.4);
\draw (9.8,2) rectangle (10.6,3.6);
\draw[dashed] (9.8,2.4) -- (10.6,2.4);
\draw[dashed] (9.8,2.8) -- (10.6,2.8);
\draw[dashed] (9.8,3.2) -- (10.6,3.2);
\draw (9,2) rectangle (11.4,2);
\draw[pattern=north west lines, pattern color=black!60] (9,4) rectangle (9.8,4.8);
\draw (9,4) rectangle (9.8,4.8);
\draw[dashed] (9,4.4) -- (9.8,4.4);
\fill[black!15] (9.8,4) -- (9.8,4.4) -- (10.6,4.4) -- (10.6,4) -- cycle;
\draw[pattern=north west lines, pattern color=black!60] (9.8,4) rectangle (10.6,4.4);
\draw (9.8,4) rectangle (10.6,4.4);
\fill[black!40] (10.6,4) -- (10.6,4.8) -- (11.4,4.8) -- (11.4,4) -- cycle;
\draw (10.6,4) rectangle (11.4,4.8);
\draw[dashed] (10.6,4.4) -- (11.4,4.4);
\draw (9,4) rectangle (11.4,4);
\draw (10.2,-0.8) node {$t=3=\tau_2\in\mathscr{T}_{H}$};
\draw (10.2,-1.3) node {$(d_3=0)$};
\draw (11,4.6) node {1};
\draw (11,4.2) node {2};
\draw (10.2,4.2) node {3};
\draw (9.4,4.6) node {4};
\draw (9.4,4.2) node {5};
\draw (10.2,1) node {1};
\draw (10.2,0.6) node {2};
\draw (10.2,0.2) node {3};
\draw (9.4,1) node {4};
\draw (9.4,0.6) node {5};
\draw (9.4,0.2) node {6};

\fill[black!15] (13.3,0) -- (13.3,1.2) -- (14.1,1.2) -- (14.1,0) -- cycle;
\draw (13.3,0) rectangle (14.1,1.2);
\draw[dashed] (13.3,0.4) -- (14.1,0.4);
\draw[dashed] (13.3,0.8) -- (14.1,0.8);
\fill[black!40] (14.1,0) -- (14.1,1.2) -- (14.9,1.2) -- (14.9,0) -- cycle;
\draw (14.1,0) rectangle (14.9,1.2);
\draw[dashed] (14.1,0.4) -- (14.9,0.4);
\draw[dashed] (14.1,0.8) -- (14.9,0.8);
\draw (12.5,0) rectangle (14.9,0);
\fill[black!15] (13.3,2) -- (13.3,3.2) -- (14.1,3.2) -- (14.1,2) -- cycle;
\draw (13.3,2) rectangle (14.1,3.2);
\draw[dashed] (13.3,2.4) -- (14.1,2.4);
\draw[dashed] (13.3,2.8) -- (14.1,2.8);
\fill[black!40] (14.1,2) -- (14.1,3.2) -- (14.9,3.2) -- (14.9,2) -- cycle;
\draw (14.1,2) rectangle (14.9,3.2);
\draw[dashed] (14.1,2.4) -- (14.9,2.4);
\draw[dashed] (14.1,2.8) -- (14.9,2.8);
\draw (12.5,2) rectangle (14.9,2);
\draw (12.5,4) rectangle (13.3,5.2);
\draw[dashed] (12.5,4.4) -- (13.3,4.4);
\draw[dashed] (12.5,4.8) -- (13.3,4.8);
\fill[black!15] (13.3,4) -- (13.3,4.8) -- (14.1,4.8) -- (14.1,4) -- cycle;
\draw[pattern=north west lines, pattern color=black!60] (13.3,4) rectangle (14.1,4.8);
\draw (13.3,4) rectangle (14.1,4.8);
\draw[dashed] (13.3,4.4) -- (14.1,4.4);
\fill[black!40] (14.1,4) -- (14.1,4.4) -- (14.9,4.4) -- (14.9,4) -- cycle;
\draw[pattern=north west lines, pattern color=black!60] (14.1,4) rectangle (14.9,4.4);
\draw (14.1,4) rectangle (14.9,4.4);
\draw (12.5,4) rectangle (14.9,4);
\draw (13.7,-0.8) node {$t=4\in\mathscr{T}_P$};
\draw (13.7,-1.3) node {$(d_4=2)$};
\draw (14.5,4.2) node {3};
\draw (13.7,4.6) node {4};
\draw (13.7,4.2) node {5};
\draw (14.5,1) node {1};
\draw (14.5,0.6) node {2};
\draw (14.5,0.2) node {3};
\draw (13.7,1) node {4};
\draw (13.7,0.6) node {5};
\draw (13.7,0.2) node {6};

\draw (2,-2.8) rectangle (2.8,-2.4);
\draw (4.75,-2.63) node {units of age 0};
\fill[black!15] (8.5,-2.8) -- (9.3,-2.8) -- (9.3,-2.4) -- (8.5,-2.4) -- cycle;
\draw (8.5,-2.8) rectangle (9.3,-2.4);
\draw (11.25,-2.63) node {units of age 1};
\fill[black!40] (2,-3.8) -- (2.8,-3.8) -- (2.8,-3.4) -- (2,-3.4) -- cycle;
\draw (2,-3.8) rectangle (2.8,-3.4);
\draw (4.75,-3.63) node {units of age 2};
\draw[pattern=north west lines, pattern color=black!60] (8.5,-3.8) rectangle (9.3,-3.4);
\draw (8.5,-3.8) rectangle (9.3,-3.4);
\draw (12.52,-3.63) node {matched units under policy $B$};
\end{tikzpicture}
\caption{An illustrative example to show the dynamic unit-matching scheme}
\label{figure:Figure2}
\end{figure}
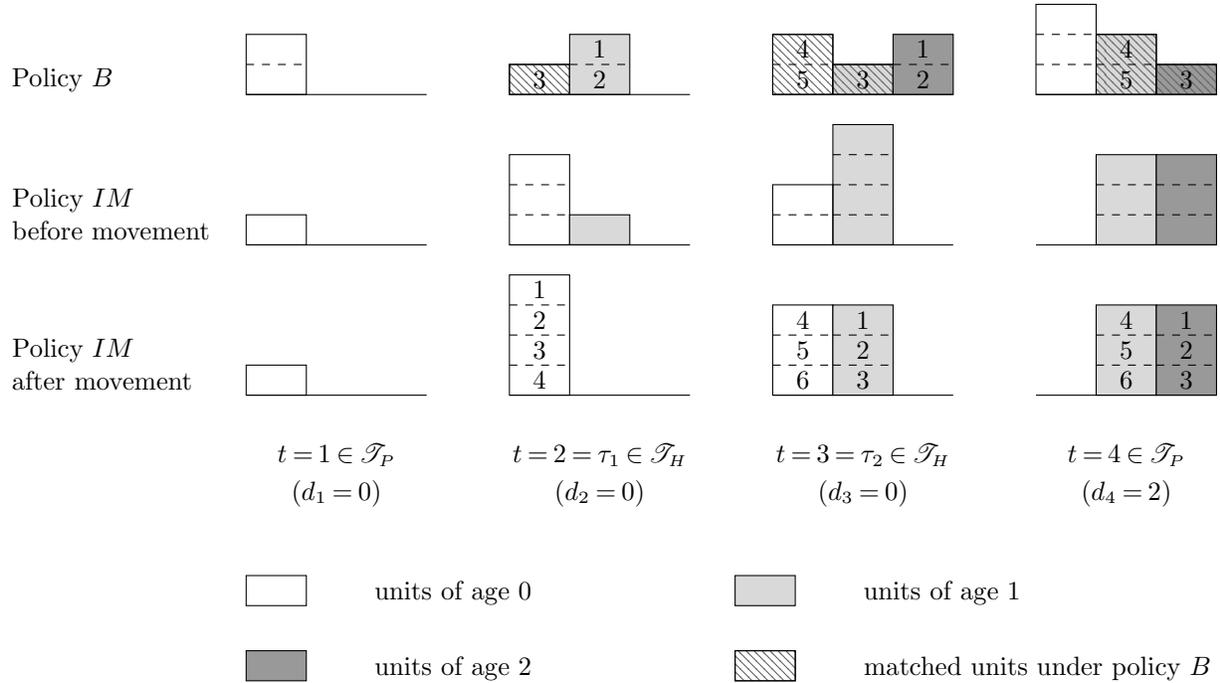

Next, we illustrate the dynamic unit-matching scheme in Figure \ref{figure:Figure2} based on the example presented in $\S$\ref{sec:im}. At the beginning of period $t=1\in \mathscr{T}_P$, no matching is defined. At the beginning of period $t=2=\tau_1$, after the movements under policy $IM$, units under policy $B$ are assigned with indices from 1 to 3, units under policy $IM$ are assigned with indices from 1 to 4, and unit 3 under policy $B$ (ordered at $\tau_1$) is temporarily matched to unit 3 under policy $IM$; {units 1-2 are not matched since they are not ordered at periods in $\mathscr{T}_H$}. At the beginning of period $t=3=\tau_2$, no permanent matching is defined since unit 3 under policy $B$ is still in inventory. After the movements under policy $IM$, units under policy $B$ are assigned with indices from 1 to 5, units under policy $IM$ are assigned with indices from 1 to 6, and units 3-5 under policy $B$ (ordered at either $\tau_1$ or $\tau_2$) are temporarily matched to units 3-5 under policy $IM$, respectively. At the beginning of period $t=4\in \mathscr{T}_P$, no new temporary matching is defined; also, since units 3-5 are all in inventory, no permanent matching is defined. At the end of the horizon (i.e., end of period 4), units 3-4 under policy $B$ is permanently matched to units 1-2 under policy $IM$ since they are used to satisfy the same units of demand. Unit 5 under policy $B$ is still in inventory at the end of the horizon, therefore it is permanently matched to its last temporarily matched unit, i.e., unit 5 under policy $IM$.

%%%%%%%%%%%%%%%%%%%%%%%%%%%%%%%%%%%%%%%%%%%%%%%%%%%%%
\subsection{Worst-Case Performance Guarantee.}
\label{sec:worst}

We now prove the worst-case performance guarantee for policy $B$. As discussed before, we use policy $IM$ as a bridging policy, and show that the expected total cost of policy $IM$ is no more that that of policy $OPT$, and the expected total cost of policy $B$ is at most twice that of policy $IM$, respectively. 

%%%%%%%%%%%%%%%%%%%%%%%%%%%%%%%%%%%%%%%%%%%%%%%%%%%%%
\textbf{{Comparison of Policies $IM$ and $OPT$}}: Recall that under Policy $IM$, an optimal ordering policy is implemented and old products can be replaced with new ones for free. Therefore, to show that the expected total cost of policy $IM$ is no more than that of policy $OPT$, it is sufficient to show that replacing old products with new ones does not increase the expected total cost. This is very intuitive {and expected in most perishable inventory problems,} because a major concern for perishable inventory systems is outdating; {and} clearly, replacing old products with new ones reduces the chance of products being outdated.

{On the other hand, while it is intuitive that replacing old products with new ones does not increase the expected total cost; it is not obvious when this would be {ensured to be} true {theoretically}. {In that regard, it turns out that this condition coincides with the optimality of FIFO issuing policy: FIFO being an optimal issuing policy implies that younger products are more preferred to have in inventory than older ones, thus replacing old products with new ones will not increase the expected total cost; and vice versa. We next formally describe the optimality of an issuing policy, followed by the statement of our assumption. 

At each period $t$, given the inventory vector $\textbf{x}_t$, the information set $f_t$ and the ordering quantity $q_t$, let $v_{k,t}$ be the amount of units of age $k$ that is used to meet demand at period $t, k=0,...,K-1$;let $x_{0,t}=q_t$, and $y_t=\sum\limits_{k=1}^{K-1}x_{k,t}+q_t$. Then, $v_{k,t}\leq x_{k,t}, k=0,...,K-1$, and $\sum\limits_{k=0}^{K-1}v_{k,t}=\min\{y_t, d_t\}$. An issuing \textit{decision rule} is a function from the set of all possible $(\textbf{x}_t,f_t,q_t,d_t)$ to the set of all possible $\textbf{v}_t=(v_{0,t},...,v_{K-1,t})$; and an issuing \textit{policy} is a collection of issuing decision rules at all periods. Then, a FIFO issuing policy is such that at each period $t$, $v_{k,t}=\min\{x_{k,t},(d_t-\sum\limits_{m=k+1}^{K-1}x_{m,t})^+\}, k=0,...,K-1$. Given an initial inventory level and an ordering policy, an issuing policy is said to be optimal if it minimizes the expected total cost among all issuing policies. Let $\Phi_t$ be the cumulative distribution function (c.d.f.) of the demand at period $t$ conditioned on $f_t$. Define the inverse of the c.d.f. $\Phi_t^{-1}(z):= \inf\{x: \Phi_t(x)\geq z\}$, and define the critical fractile $\bar{y}_t:=\Phi_t^{-1}(\frac{p}{p+h})$. To compare policies $IM$ and $OPT$, we state our assumption as follows:

\begin{assumption}
\label{assump:1}
Starting from any period $t$, suppose $y_t\leq \max\limits_{\tau=1,...,t}\bar{y}_{\tau}$ holds at $t$ and an optimal ordering policy is implemented at $t+1,...,T$, then FIFO is an optimal issuing policy, i.e., FIFO minimizes the expected total cost at $t,...,T$ among all issuing policies. 
\end{assumption}

{Intuitively, this assumption says that if the total inventory level after ordering is less than or equal to a threshold at some period, and an optimal ordering policy is implemented at the following periods, then FIFO issuing policy minimizes the expected total cost from that period to the end of the planning horizon.} Note that a stronger sense of optimality of an issuing policy requires the issuing policy to be optimal for any initial inventory level and under an arbitrary sequence of ordering quantities (\citet{pierskalla1972optimal}). However, we only require FIFO to be optimal for small initial inventory levels and under an optimal ordering policy, which is a much weaker assumption.

We acknowledge that in practice, directly checking whether Assumption \ref{assump:1} holds or not may be difficult. To date, FIFO has already been shown to be optimal under i.i.d. demand (\citet{fries1975optimal2}), which is widely assumed in the perishable inventory literature, or zero holding cost (\citet{pierskalla1972optimal}). In $\S$\ref{sec:fifo}, we extend these existing findings and present a necessary and sufficient condition and three easy-to-check sufficient conditions to ensure the optimality of FIFO issuing policy. In particular, we show that the condition presented in Chao et al. (2015) to ensure a performance guarantee of two (i.e., demand is independent and stochastically non-decreasing over time) is a special case of ours, and we further provide conditions and examples where FIFO is optimal and our performance guarantee is strictly tighter.

With Assumption \ref{assump:1}, we now present a structural property on the optimal cost-to-go function, which is a key result for comparing policies $IM$ and $OPT$. At each period $t$, given $\textbf{x}_t$ and $f_t$, let $C_t(\textbf{x}_t,f_t)$ be the optimal cost-to-go, and let $y_t=\sum\limits_{k=1}^{K-1}x_{k,t}+q_t$. Then, we have:
$$C_t(\textbf{x}_t, f_t)=\min_{q_t\geq 0}\{p\mathrm{E}[(D_t-y_t)^+|f_t]+h\mathrm{E}[(y_t-D_t)^+|f_t]+w\mathrm{E}[(x_{K-1,t}-D_t)^+|f_t]$$
\begin{equation*}
+\beta \mathrm{E}[C_{t+1}(\textbf{X}_{t+1},F_{t+1})]\}.
\end{equation*}
For $k=1,...,K-1$, for the continuous case, let $C_t^{(k)}(\textbf{x}_t,f_t)$ denote the partial derivative of $C_t(\textbf{x}_t,f_t)$ with respect to to $x_{k,t}$ (the differentiability can be easily established following similar arguments as in \citet{fries1975optimal2}); for the discrete case, let $C_t^{(k)}(\textbf{x}_t,f_t)$ denote the incremental of $C_t(\textbf{x}_t,f_t)$ caused by a unit increase of $x_{k,t}$. Then, we have the following result.

\begin{lemma}
\label{lem:ij}
Under Assumption \ref{assump:1}, for $t=1,...,T$, (i) $C_{t+1}^{(k)}(\textbf{x}_{t+1},f_{t+1})\geq 0, k=1,...,K-1$, $\forall \textbf{x}_{t+1},f_{t+1}$; (ii) $C_{t+1}^{(i)}(\textbf{x}_{t+1},f_{t+1})\leq C_{t+1}^{(j)}(\textbf{x}_{t+1},f_{t+1})\leq w/\beta, 1\leq i<j\leq K-1$, $\forall \textbf{x}_{t+1},f_{t+1}$ such that $\sum\limits_{k=1}^{K-1}x_{k,t+1}< \max\limits_{\tau=1,...,t}\bar{y}_{\tau}-d_{t}$.
\end{lemma}

Lemma \ref{lem:ij} implies that when FIFO is an optimal issuing policy (i.e., Assumption \ref{assump:1}), the optimal cost-to-go is non-decreasing in the inventory levels; and if the total inventory level is small, the incremental of the discounted optimal cost-to-go caused by a unit increase of an older unit is higher than that caused by a younger unit, and both are bounded by the unit outdating cost. Then, since units are only moved from older to younger positions under policy $IM$ {(it is not difficult to show that $\sum\limits_{k=1}^{K-1}x_{k,t+1}^{IM}\leq (\max\limits_{\tau=1,...,t}\bar{y}_{\tau}-d_{t})^+$ holds for all $t$)}, it is intuitive that the expected total cost of policy $IM$ is no more than that of policy $OPT$, which we formally prove in the next lemma.

\begin{lemma}
\label{lem:optim}
Under Assumption \ref{assump:1}, $E[\mathscr{C}(IM)]\leq E[\mathscr{C}(OPT)]$.
\end{lemma}

%%%%%%%%%%%%%%%%%%%%%%%%%%%%%%%%%%%%%%%%%%%%%%%%%%%%%
\textbf{Comparison of Policies $B$ and $IM$}: With Lemma \ref{lem:optim}, to establish the worst-case performance guarantee of policy $B$, all remains to show is $\mathrm{E}[\mathscr{C}(B)]\leq 2\mathrm{E}[\mathscr{C}(IM)]$. To do so, we first provide a key lemma as follows, where $P_t^{\pi}, H_t^{\pi}$ and $W_t^{\pi}$ denote the marginal shortage penalty, holding and outdating costs at period $t$ under policy $\pi$, respectively.

\begin{lemma}
\label{lem:imb1}
With probability one, (i) $\sum\limits_{t\in \mathscr{T}_P}P_t^{B}\leq \sum\limits_{t=1}^TP_t^{IM}$; (ii) $\sum\limits_{t\in \mathscr{T}_H}H_t^{B}\leq \sum\limits_{t=1}^TH_t^{IM}$; (iii) $\sum\limits_{t\in \mathscr{T}_H}W_t^{B}\leq \sum\limits_{t=1}^TW_t^{IM}.$
\end{lemma}

\proof{Proof.} (i) For any given sample path, after the movements of units under policy $IM$, we have $y_t^B\geq y_t^{IM}, \forall t\in \mathscr{T}_P$. Then, $P_t^B\leq P_t^{IM}, \forall t\in \mathscr{T}_P$. Therefore, $\sum\limits_{t\in \mathscr{T}_P}P_t^{B}\leq \sum\limits_{t\in \mathscr{T}_P}P_t^{IM}\leq \sum\limits_{t=1}^TP_t^{IM}$ with probability one.

(ii) To show that the total holding cost charged for units ordered at $t\in\mathscr{T}_H$ under policy $B$ is no more than the total holding cost under policy $IM$, it is sufficient to show that under the dynamic unit-matching scheme described in $\S$\ref{sec:matching}, all matched units under policy $B$ (i.e., all units ordered at $t\in\mathscr{T}_H$) stay in inventory no longer than their permanently matched units under policy $IM$. 

Units ordered at $t\in\mathscr{T}_H$ under policy $B$ can be 1) used to satisfy demand, 2) outdated or still in inventory at the end of the planning horizon. First, recall that after the movements at each period $\tau_i\in \mathscr{T}_H$, units under policy $B$ are assigned indices from 1 to $y_{\tau_i}^{B}$ and units under policy $IM$ are assigned indices from 1 to $y_{\tau_i}^{IM}$ from oldest to youngest. Since the inventory vector under policy $IM$ is ``younger'' than that under policy $B$ (i.e., Inequality (\ref{Eq. ine})), each unit under policy $IM$ is consumed no earlier than the unit with the same index under policy $B$. Therefore, for a matched unit under policy $B$ that is used to satisfy demand (call it $u_1$), its last temporarily matched unit under policy $IM$, which has the same index and age as $u_1$, is consumed no earlier than $u_1$. Then, given the FIFO issuing policy, the permanently matched unit of $u_1$, which is used to satisfy the same unit of demand as $u_1$, must be no younger than $u_1$.

Second, by the dynamic unit-matching scheme, for a matched unit under policy $B$ that is outdated or is still in inventory at the end of the planning horizon (call it $u_2$), its permanently matched unit under policy $IM$ is defined as its last temporarily matched unit. Since all temporarily matched pairs of units have the same age, and further considering possible movements from older to younger positions for units under policy $IM$, $u_2$ stays in inventory no longer than its permanently matched unit. 

Since the permanent matchings are defined on a one-to-one correspondence and the above argument is true for any given sample path, we have $\sum\limits_{t\in \mathscr{T}_H}H_t^{B}\leq \sum\limits_{t=1}^TH_t^{IM}$ with probability one.

(iii) By the dynamic unit-matching scheme, for a matched unit under policy $B$ that is outdated, the permanently matched unit under policy $IM$ must be outdated at the same period. Since the permanent matchings are defined on a one-to-one correspondence and the above argument is true for any given sample path, we have $\sum\limits_{t\in \mathscr{T}_H}W_t^{B}\leq \sum\limits_{t=1}^TW_t^{IM}$ with probability one.
\Halmos
\endproof

With the above result, now it is easy to reach to the following conclusion.

\begin{lemma}
\label{lem:imb}
$\mathrm{E}[\mathscr{C}(B)]\leq 2\mathrm{E}[\mathscr{C}(IM)]$.
\end{lemma}

\proof{Proof.} Let $\mathbbm{1}{(t\in \mathscr{T}_P)}$ and $\mathbbm{1}{(t \in \mathscr{T}_H)}$ be two indicator functions. Then, we have $\mathbbm{1}{(t\in \mathscr{T}_P)}+\mathbbm{1}{(t \in \mathscr{T}_H)}=1$ with probability one, and we have the following result.
\begin{align*}
\mathrm{E}[\mathscr{C}(B)]
=&\sum\limits_{t=1}^T\mathrm{E}[\mathrm{E}[(P_t^{B}+H_t^{B}+W_t^{B})|F_t]]\\
=&\sum\limits_{t=1}^T\mathrm{E}[\mathrm{E}[(P_t^{B}+H_t^{B}+W_t^{B})(\mathbbm{1}{(t\in \mathscr{T}_P})+\mathbbm{1}{(t\in \mathscr{T}_H)})|F_t]] \notag \\
=&\sum\limits_{t=1}^T\mathrm{E}[\mathrm{E}[2P_t^B\mathbbm{1}{(t\in \mathscr{T}_P})+2(W_t^B+H_t^B)\mathbbm{1}{(t\in \mathscr{T}_H)}|F_t]] \notag \\
=&\mathrm{E}\bigg[\sum\limits_{t\in \mathscr{T}_P}2P_t^{B}+\sum\limits_{t\in \mathscr{T}_H}2(H_t^{B}+W_t^{B})\bigg]\\
\leq &\mathrm{E}\bigg[2\sum\limits_{t=1}^TP_t^{IM}+2\sum\limits_{t=1}^T(H_t^{IM}+W_t^{IM})\bigg]\\
=&2\mathrm{E}[\mathscr{C}(IM)],
\end{align*}
where the third equality follows from the definition of the marginal-cost dual-balancing policy and the fact that $\mathbbm{1}{(t\in \mathscr{T}_P)}$ and $\mathbbm{1}{(t \in \mathscr{T}_H}$ are deterministic for given $f_t$, and the inequality follows from Lemma \ref{lem:imb1}.
\Halmos
\endproof

Based on the above results, we now state our main theorem as follows.

\begin{theorem}
\label{optb}
Under Assumption \ref{assump:1}, the marginal-cost dual-balancing policy has a worst-case performance guarantee of two. That is, the expected total cost of the marginal-cost dual-balancing policy is at most twice that of an optimal ordering policy, i.e., $\mathrm{E}[\mathscr{C}(B)]\leq 2\mathrm{E}[\mathscr{C}(OPT)]$.
\end{theorem}

\proof{Proof.}
Since we have $\mathrm{E}[\mathscr{C}(IM)]\leq \mathrm{E}[\mathscr{C}(OPT)]$ from Lemma \ref{lem:optim} and $\mathrm{E}[\mathscr{C}(B)]\leq 2\mathrm{E}[\mathscr{C}(IM)]$ from Lemma \ref{lem:imb}, we have $\mathrm{E}[\mathscr{C}(B)]\leq 2\mathrm{E}[\mathscr{C}(IM)] \leq 2\mathrm{E}[\mathscr{C}(OPT)]$, which completes the proof. \Halmos
\endproof

%%%%%%%%%%%%%%%%%%%%%%%%%%%%%%%%%%%%%%%%%%%%%%%%%%%%%
\subsection{{Comparison with an Optimal Base-Stock Policy.}}
\label{sec:base}

In the previous subsection, we have compared our policy with an optimal ordering policy, and have shown that the expected total cost of our policy is at most twice that of an optimal ordering policy when FIFO is an optimal issuing policy (i.e., Assumption \ref{assump:1}). In this subsection, we switch the benchmark to compare our policy with an optimal base-stock policy, and show that the expected total cost of our policy is always at most twice that of an optimal base-stock policy.

This result is noteworthy because as discussed in $\S$\ref{sec:intro}, among the heuristic policies developed for perishable inventory systems, the base-stock policy, under which the total inventory is replenished up to the same level at each period, is particularly popular due to its simple implementation and competitive numerical performance (e.g., \citet{nahmias1976myopic, cohen1976analysis, chazan1977markovian, nandakumar1993near, cooper2001pathwise, li2009note, chen2014coordinating, zhang2016nonparametric}). However, the computation of an optimal base-stock policy involves the evaluation of the expected outdating cost for each given base-stock level, which is again intractable due to the large state space. Although many heuristic approaches are developed to compute ``good'' base-stock levels, none of them admits a theoretical performance guarantee.

In the following theorem, we compare our marginal-cost dual-balancing policy with an optimal base-stock policy (denoted as $BA$), and show that the expected total cost of our policy {(although itself is not a base-stock policy)} is at most twice that of an optimal base-stock policy. We remark that since FIFO is always optimal under base-stock policies (\citet{chazan1977markovian}), Theorem \ref{baseb} follows without any assumption.
\begin{lemma}[Chazan and Gal (1977)]
\label{lem:base}
Under base-stock policies, the cumulative amount of outdate under FIFO issuing policy is smaller than that under any other issuing policy with probability one.
\end{lemma}
\begin{theorem}
\label{baseb}
The expected total cost of the marginal-cost dual-balancing policy is at most twice that of an optimal base-stock policy, i.e., $\mathrm{E}[\mathscr{C}(B)]\leq 2\mathrm{E}[\mathscr{C}(BA)]$.
\end{theorem}

The proof for Theorem \ref{baseb} is similar to that for Theorem \ref{optb}, except that now policy $IM$ is constructed based on $BA$ instead of $OPT$. Since FIFO is always optimal under base-stock policies (Lemma \ref{lem:base}), we can easily show that $\mathrm{E}[\mathscr{C}(IM)]\leq \mathrm{E}[\mathscr{C}(BA)]$ (in fact, we have $\mathscr{C}(IM)\leq \mathscr{C}(BA)$ with probability one). Given this result and also $\mathrm{E}[\mathscr{C}(B)]\leq 2\mathrm{E}[\mathscr{C}(IM)]$ as we have shown in Lemma \ref{lem:imb}, the result in Theorem \ref{baseb} immediately follows.

Also note that the result in Theorem $\ref{baseb}$ can be easily extended to cases where the base-stock levels at different periods are different but non-decreasing over time.

%%%%%%%%%%%%%%%%%%%%%%%%%%%%%%%%%%%%%%%%%%%%%%%%%%%%
\section{Truncated-Balancing Policy.}
\label{sec:truncated}

{In this section, we first prove that a myopic policy under the marginal-cost accounting scheme provides a lower bound on the optimal ordering quantity. Then, by combining the \textit{specific} lower bound we derive and any upper bound on the optimal ordering quantity with the marginal-cost dual-balancing policy, we present a more general class of algorithms that we call the truncated-balancing policy.} We later show that while both the marginal-cost dual-balancing policy and the truncated-balancing policy have the same worst-case performance guarantee of two, the latter performs much better in the computational studies (see $\S$\ref{sec:numerical}).

In the following proposition, we show that under Assumption \ref{assump:1}, the minimizer of the total marginal costs provides a lower bound on the optimal ordering quantity. A similar result has been developed in \citet{levi2007approximation} for the nonperishable backlogging case (clearly, Assumption \ref{assump:1} is not needed for the nonperishable case), where the total inventory level is known to be a sufficient statistic for the system state. However, generalizing this result to the lost sales case, where a pipeline inventory vector is needed to describe the system state, remains an open problem. In that regard, our problem is similar to the lost sales case since we also need an inventory vector to describe the system state, and the analysis for the nonperishable backlogging case is not applicable to our case. 

\begin{proposition}
\label{prop:lower}
At any period $t$, given $\textbf{x}_t$ and $f_t$, let $q_t^{OPT}$ be an optimal ordering quantity, and $q_t^{L}$ be the smallest quantity that minimizes $P_t(\textbf{x}_t,f_t,q_t)+H_t(\textbf{x}_t,f_t,q_t)+W_t(\textbf{x}_t,f_t,q_t)$. Then, under Assumption \ref{assump:1}, $q_t^{L}\leq q_t^{OPT}$.
\end{proposition}

\begin{remark}
The reason why a myopic ordering quantity under the marginal-cost accounting scheme provides a lower bound on the optimal ordering quantity is that, under the marginal-cost accounting scheme, the optimal cost-to-go function is non-increasing in the inventory levels of all ages (see proof in the Appendix); thus, ordering more units can decrease the optimal cost-to-go of the next period, and the minimizer of the single-period cost, which ignores the benefit of ordering more for future, will tend to order less than optimal. We remark that while the monotonicity result for the optimal cost-to-go function can be easily established for the nonperishable inventory case, it can in fact be violated for the perishable inventory case in general. However, as we show in the proof of Proposition \ref{prop:lower}, the optimality of FIFO issuing policy is a sufficient condition to establish the monotonicity result, which we believe is a new and important contribution to the literature.
\end{remark}

Based on the above result, we now define the \textit{truncated-balancing policy} (denoted as $TB$) as follows. At each period $t$, given $\textbf{x}_t$ and $f_t$, let $q_t^B$ be the balancing ordering quantity defined by Equation (\ref{eqn:balancing}), let $q_t^L$ be the lower bound on the optimal ordering quantity defined in Proposition \ref{prop:lower} (or any \textit{looser} lower bound), and let $q_t^U$ be any upper bound on the optimal ordering quantity. Then, the truncated-balancing ordering quantity $q_t^{TB}$ is defined as:
$$q_t^{TB}=
\left\{
\begin{aligned}
&q_t^B,\quad \text{if}~q_t^L\leq q_t^B\leq q_t^U,\\
&q_t^L,\quad \text{if}~q_t^B<q_t^L,\\
&q_t^U,\quad \text{if}~q_t^B>q_t^U.
\end{aligned}
\right.
$$

Note that for trivial lower and upper bounds (i.e., $q_t^L=0, q_t^U=\infty$), the truncated-balancing policy reduces to the marginal-cost dual-balancing policy. In the following theorem, we show that the truncated-balancing policy, as a more general class of algorithms, also admits a worst-case performance guarantee of two.

\begin{theorem}
\label{opttb}
Under Assumption \ref{assump:1}, the truncated-balancing policy has a worst-case performance guarantee of two. That is, the expected total cost of the truncated-balancing policy is at most twice that of an optimal ordering policy, i.e., $\mathrm{E}[\mathscr{C}(TB)]\leq 2\mathrm{E}[\mathscr{C}(OPT)]$.
\end{theorem}

\begin{remark}
{We remark that policy $TB$ is not guaranteed to perform at least as good as policy $B$ (thus the proof of Theorem \ref{opttb} is nontrivial). While it may appear that $q_t^{TB}$ is at least as good as $q_t^B$, this is only true if an optimal policy is implemented at the following periods.}
\end{remark}

\begin{remark}{We also remark that unlike the nonperishable inventory case where the lower and upper bounds in the definition of policy $TB$ can be replaced with any (tighter) ones, in our case, the lower bound $q_t^L$, as a minimizer of the single-period marginal cost, is special and cannot be tightened. To see why this is the case, let $y_t^{TB}$ and $y_t^{IM}$ be the total inventory levels after ordering at period $t$ under policies $TB$ and $IM$, respectively; also, given $\textbf{x}_t^{TB}$ and $f_t$, let $y_t^B$ denote the total inventory level after ordering if the balancing ordering quantity $q_t^B$ is ordered. Then, given that $q_t^B<q_t^{TB}=q_t^L$, it is possible that $y_t^B\leq y_t^{IM}<y_t^{TB}$ (while for the nonperishable case, since a base-stock policy is optimal, given that $q_t^B<q_t^{TB}=q_t^L$, we always have $y_t^B< y_t^{TB}\leq y_t^{OPT}$). In this case, since $y_t^{TB}>y_t^{IM}$, it is not possible to match all the $q_t^{TB}$ units to units under policy $IM$. Therefore, we instead only match the first $q_t^B$ units ordered under policy $TB$ to units under policy $IM$, and then show that the total marginal cost at period $t$ for ordering $q_t^{TB}$ is no more than that for ordering $q_t^B$, which is only ensured to be true if $q_t^{TB}=q_t^L$ is a minimizer of the single-period marginal cost (see more details in the proof of Theorem \ref{opttb} in the Appendix).}
\end{remark}

%%%%%%%%%%%%%%%%%%%%%%%%%%%%%%%%%%%%%%%%%%%%%%%%%%%%%%
\section{{Sufficient Conditions for Optimality of FIFO Issuing Policy.}}
\label{sec:fifo}

In $\S$\ref{sec:analysis}-\ref{sec:truncated}, we have shown that our proposed algorithms have a worst-case performance guarantee of two when FIFO is an optimal issuing policy (i.e., Assumption \ref{assump:1}). In this section, we provide a necessary and sufficient condition and several easy-to-check sufficient conditions that ensure the optimality of FIFO issuing policy, which extends the existing literature and provides insights into the key trade-offs of different issuing policies. 

In Lemma \ref{lem:ij}, we have presented a necessary condition for the optimality of FIFO issuing policy. In the following proposition, we show that this condition is also sufficient (in fact, we only need a part of that condition), which leads to a \textit{necessary and sufficient} condition for the optimality of FIFO issuing policy.
\begin{proposition}
\label{prop:fifo}
Assumption \ref{assump:1} holds if and only if for $t=1,...,T$, $C_{t+1}^{(k)}(\textbf{x}_{t+1},f_{t+1})\leq w/\beta, k=1,...,K-1$, $\forall \textbf{x}_{t+1},f_{t+1}$ such that $\sum\limits_{k=1}^{K-1}x_{k,t+1}< \max\limits_{\tau=1,...,t}\bar{y}_{\tau}-d_{t}.$
\end{proposition}

Proposition \ref{prop:fifo} says that FIFO is optimal if and only if the incremental of the discounted optimal cost-to-go caused by a unit increase of the inventory of any age is bounded by the unit outdating cost. This provides an overall insight into the key trade-off in ensuring the optimality of FIFO issuing policy, and provides intuition in finding sufficient conditions for Assumption \ref{assump:1} to hold. In particular, a unit increase of the in-hand inventory can potentially increase both holding and outdating costs but decrease the shortage penalty for future periods. Consider a case where demand for future periods is sufficiently large so that the decrease in shortage penalty offsets the increase in holding cost. In this case, the incremental of the total cost is bounded by the unit outdating cost, i.e., the condition in Proposition \ref{prop:fifo} holds. Also, consider another case where the unit holding cost is sufficiently small so that the total discounted holding and outdating costs for future periods are bounded by the unit outdating cost, regardless of the demand distribution. In this case, we also have the condition in Proposition \ref{prop:fifo} hold.

Based on these interpretations of Proposition \ref{prop:fifo}, we now provide several easy-to-check sufficient conditions that ensure the optimality of FIFO issuing policy. In particular, in Propositions \ref{prop:demand} and \ref{prop:cost}, we formalize the two intuitions discussed above and show that Assumption \ref{assump:1} holds either when future demand is sufficiently large or when holding cost is sufficiently small. Next, in Proposition \ref{prop:mix}, we present a sufficient condition for the optimality of FIFO issuing policy that involves both demand and holding cost; however the respective conditions are much weaker than those that involve either only demand (Proposition \ref{prop:demand}) or only holding cost (Proposition \ref{prop:cost}). {We remark that our results also extend the existing findings on the optimality of FIFO issuing policy (\citet{fries1975optimal2, pierskalla1972optimal}), which we describe in further details below.}

In the following Proposition, we first present a sufficient condition on demand distribution that ensures the optimality of FIFO issuing policy.

\begin{proposition}
\label{prop:demand}
Assumption \ref{assump:1} holds if $\forall f_{T+1}, \bar{y}_t=\Phi_{t}^{-1}(\frac{p}{p+h})$ is non-decreasing in $t$.
\end{proposition}

As discussed above, the intuition for the above condition is that when future demand is large, increasing the in-hand inventory level will not increase the sum of shortage penalty and holding cost, and thus the condition in Proposition \ref{prop:fifo} holds, which implies the optimality of FIFO issuing policy. Clearly, the condition in Proposition \ref{prop:demand} is an extension of i.i.d. demand considered in \citet{fries1975optimal2}, and also includes the independent and stochastically non-decreasing demand (in the sense of first-order dominance)\footnote{Note that the exact condition \citet{chao2015approximation} needs in their proof is that $S_t$ is non-decreasing in $t$, where $S_t$ is the solution of y to the equation $\mathrm{E}[(y-D_t)^+|f_t]=\mathrm{E}[(D_t-y)^+|f_t]$.} presented in \citet{chao2015approximation} as a special case. 

Next, we present a sufficient condition regarding holding cost to ensure the optimality of FIFO issuing policy, which extends the zero holding cost considered in \citet{pierskalla1972optimal}.

\begin{proposition}
\label{prop:cost}
Assumption \ref{assump:1} holds if $h\leq \frac{1-\beta}{\beta}w$.
\end{proposition}

This result is also intuitive because when $h\leq \frac{(1-\beta)}{\beta}w$, the total discounted holding cost for an arbitrary number of future periods $(\beta+\beta^2+...)h=\beta h/(1-\beta)$ is bounded by the unit outdating cost $w$. 

Note that Proposition \ref{prop:cost} identifies an important class of problems where the performance guarantee of our algorithms is strictly tighter than that presented in \citet{chao2015approximation}. Recall the cost transformation that $h=\hat{h}+(1-\beta)\hat{c}, w=\hat{w}+\beta \hat{c}$. Then it is straightforward to check that $h\leq \frac{1-\beta}{\beta}w$ if and only if $\hat{h}\leq \frac{1-\beta}{\beta}\hat{w}$. Thus, Proposition \ref{prop:cost} implies that as long as the original holding cost $\hat{h}$ is sufficiently small, our algorithms have a performance guarantee of two. Under general demand, the performance guarantee presented in \citet{chao2015approximation} is $2+\frac{(K-2)h}{Kh+w}$. For general $K>2$ and any $\beta\in(0,1)$, for their performance guarantee $2+\frac{(K-2)h}{Kh+w}=2+\frac{(K-2)(\hat{h}+(1-\beta)\hat{c})}{K(\hat{h}+(1-\beta)\hat{c})+\hat{w}+\beta \hat{c}}$ to be equal to 2, the transformed holding cost $h$ needs to be zero, which implies both the original ordering cost $\hat{c}$ and holding cost $\hat{h}$ need to be zero. On the other hand, we only assume the original holding cost $\hat{h}$ to be small and allow arbitrarily large ordering cost $\hat{c}$ (note that conducting cost transformation to eliminate $\hat{c}$ with end up with positive $h$). Therefore, {as we illustrate in the following example,} the performance guarantee of our algorithms can be strictly tighter, especially for problems with positive ordering cost and small or zero holding cost, which represent a large class of perishable (especially blood) inventory problems (\citet{haijema2007blood,zhou2011inventory}). 

\begin{example}
\label{ex:1}
Consider an instance with product lifetime $K=5$, and original cost parameters $\hat{c}=1, \hat{h}=\hat{w}=0$. Demand is a general stochastic process. Then for discount factor values $\beta=0.9,0.95,0.99$, the performance guarantees presented in \citet{chao2015approximation} are 2.214, 2.125, 2.029, respectively, while the performance guarantee of our algorithms is exactly two in all three cases.
\end{example}

Next, inspired by Propositions \ref{prop:demand} and \ref{prop:cost}, in the following proposition, we present a sufficient condition for the optimality of FIFO issuing policy that involves both demand and holding cost. We further show that Proposition \ref{prop:mix} is in fact much more general and provides a unified framework that extends both Propositions \ref{prop:demand} and \ref{prop:cost}.

\begin{proposition}
\label{prop:mix}
Assumption \ref{assump:1} holds if $\forall f_{T+1}, h\leq \frac{1-\gamma}{\gamma}p+\frac{1-\beta \gamma}{\beta \gamma}w$, where $\gamma=\max\limits_{1< s\leq t\leq T}\Phi_t(\bar{y}_s)$.
\end{proposition}

Clearly, $\gamma$ provides an upper bound on the probability that there will be excess inventory after demand realization at each period. The intuition behind Proposition \ref{prop:mix} is that by increasing a unit of in-hand inventory and keeping the ordering quantity unchanged, the increase of the total expected holding and outdating costs are at most $\gamma(h+w)$ while the decrease of the expected shortage penalty is at least $(1-\gamma)p$. Then, the incremental of the expected total cost by increasing a unit of in-hand inventory is bounded by $\gamma(h+w)-(1-\gamma)p$; let it to be less than or equal to $w/\beta$ (in order to satisfy the condition in Proposition \ref{prop:fifo}) and rearrange terms, we achieve the condition in Proposition \ref{prop:mix}. 

We believe this is an authentic result and provides key insights into the main trade-offs of FIFO issuing policy. In particular, Proposition \ref{prop:mix} says that FIFO is optimal when the demand over time does not ``drop'' significantly (so that $\gamma$ is not too large) and the holding cost is moderately small. Moreover, there is a clear trade-off between demand and holding cost: The smaller the holding cost, the fewer requirements we need for the demand, and vice versa. This result is powerful and in fact provides a unified framework that extends both Propositions \ref{prop:demand} and \ref{prop:cost}:} On one hand, suppose $\forall f_{T+1}$, the critical fractile $\bar{y}_t=\Phi_{t}^{-1}(\frac{p}{p+h})$ is non-decreasing in $t$; then $\gamma=\Phi_t(\bar{y}_t)=\frac{p}{p+h}$ and thus $h(= \frac{1-\gamma}{\gamma}p)\leq \frac{1-\gamma}{\gamma}p+\frac{1-\beta \gamma}{\beta \gamma}w$ is automatically satisfied, i.e., Proposition \ref{prop:mix} reduces to Proposition \ref{prop:demand}. On another hand, suppose we do not impose any restriction on demand; then the worst case is $\gamma=1$ and we thus need $h\leq \frac{1-\beta}{\beta}w$, i.e., Proposition \ref{prop:mix} reduces to Proposition \ref{prop:cost}.

Furthermore, the condition in Proposition \ref{prop:mix} also leads to a broader class of problems where the performance guarantee of our algorithms is strictly tighter than that in \citet{chao2015approximation}, {as we illustrate in the following example}.

\begin{example}
\label{ex:2}
Consider an instance with product lifetime $K=5$, cost parameters $p=w=5,h=1$, and discount factor $\beta=1$ (any $\beta\leq 1$ would work). Let the demand at each period be independent and assume demands at periods $1,3,5,...$ are exponentially distributed with mean of 5, i.e., exp($\frac{1}{5}$), and demands at periods $2,4,6,...$ are exp($\frac{1}{6}$) (note that our result does not rely on this specific pattern of fluctuation; the demand can be for example exp($\frac{1}{5}$) on weekdays and exp($\frac{1}{6}$) on weekends). It is easy to check that the critical fractile for exp($\frac{1}{5}$) is $\Phi_{1}^{-1}(\frac{p}{p+h})\approx\Phi_{1}^{-1}(0.833)\approx8.9$ and the critical fractile for exp($\frac{1}{6}$) is $\Phi_{2}^{-1}(\frac{p}{p+h})\approx\Phi_{2}^{-1}(0.833)\approx10.7$. Then, $\gamma=\max\limits_{1< s\leq t\leq T}\Phi_t(\bar{y}_s)\approx\Phi_{1}(10.7)\approx0.882$, and $h=1< \frac{1-\gamma}{\gamma}p+\frac{1-\beta \gamma}{\beta \gamma}w\approx1.338$. Thus, the condition in Proposition \ref{prop:mix} is clearly satisfied and the performance guarantee of our algorithms is two. On the other hand, the performance guarantee presented in \citet{chao2015approximation} is $2+\frac{(K-2)h}{Kh+w}= 2.3$.
\end{example}

We remark that besides the sufficient conditions we presented above, there could be potentially many other situations where FIFO issuing policy is optimal (i.e., the condition in Proposition \ref{prop:fifo} holds). It is worth noting that under all those situations, our algorithms would admit a worst-case performance guarantee of two.

%%%%%%%%%%%%%%%%%%%%%%%%%%%%%%%%%%%%%%%%%%%%%%%%%%%%%%
\section{Computational Experiments.}
\label{sec:numerical}

We start with a discussion on the computation of the balancing ordering quantity in $\S$\ref{sec:computation}. Then, in $\S$\ref{sec:compound}, we test the performances of our policies under the Hospital Alpha platelet inventory control problem using real data. 

%%%%%%%%%%%%%%%%%%%%%%%%%%%%%%%%%%%%%%%%%%%%%%%%%%%%%
\subsection{Computation of Balancing Ordering Quantity.}
\label{sec:computation}

At each period $t$, given $\textbf{x}_t$ and $f_t$, the marginal shortage penalty $P_t(\textbf{x}_t,f_t,q_t)$ is non-increasing in $q_t$, while the marginal holding and outdating costs $H_t(\textbf{x}_t,f_t,q_t)$ and $W_t(\textbf{x}_t,f_t,q_t)$ are non-decreasing in $q_t$. Thus the balancing ordering quantity $q_t^B$ defined in Equation (\ref{eqn:balancing}) can be computed using a simple binary search. However, to do so, we first need to efficiently compute the expected marginal costs for each given $q_t$. Given the distribution of $D_t$, the computation of $P_t(\textbf{x}_t,f_t,q_t)$ is straightforward. Thus, in this subsection, we focus on the computation of $H_t(\textbf{x}_t,f_t,q_t)$ and $W_t(\textbf{x}_t,f_t,q_t)$.

Similar as the existing studies on balancing policies, for general demands, the expected marginal holding and outdating costs $H_t(\textbf{x}_t,f_t,q_t)$ and $W_t(\textbf{x}_t,f_t,q_t)$ can be computed using methods such as Monte Carlo simulation. However, if demand over time is independent and integer-valued (we consider integer-valued quantities in our computational experiments), we can further achieve closed-form expressions for the expected marginal costs as follows.

Recall that $A_{0,t}=0$, and for $k=1,...,K-1$, $A_{k,t}$ denotes the total demand over periods $t, ...,t+k-1$ that cannot be satisfied by the inventory of $(x_{K-k,t},...,x_{K-1,t})$, and $A_{0,t}=0$. Similar to \citet{nahmias1975optimal}, for given $(x_{K-k+1,t},...,x_{K-1,t})$ and $f_t$, define:
$$R_{k,t}(x_{K-k,t})=\mathrm{P}(A_{k-1,t}+D_{t+k-1}< x_{K-k,t}|f_t), k=1,...,K,$$
which denotes the conditional probability that there will be outdates at the end of period $t+k-1$. Then, $\forall u\geq 0$, we have $R_{1,t}(u)=\mathrm{P}(D_{t}< u|f_t)$, and for $k=2,...,K$:
\begin{align*}
R_{k,t}(u)&=\mathrm{P}(A_{k-1,t}+D_{t+k-1}< u|f_t)\\
&=\sum\limits_{v=1}^{u}\mathrm{P}(A_{k-1,t}< v|f_t)\mathrm{P}(D_{t+k-1}=u-v|f_t)\\
&=\sum\limits_{v=1}^{u}\mathrm{P}((A_{k-2,t}+D_{t+k-2}-x_{K-k+1,t})^+< v|f_t)\mathrm{P}(D_{t+k-1}=u-v|f_t)\\
&=\sum\limits_{v=1}^{u}\mathrm{P}(A_{k-2,t}+D_{t+k-2}< v+x_{K-k+1,t}|f_t)\mathrm{P}(D_{t+k-1}=u-v|f_t)\\
&=\sum\limits_{v=1}^{u}R_{k-1,t}(v+x_{K-k+1,t})\mathrm{P}(D_{t+k-1}=u-v|f_t).
\end{align*}

Therefore, the probabilities $R_{k,t}$ can be computed efficiently by recursion. In this case, at each period $t$, given $\textbf{x}_t$, $f_t$ and $q_t$, the expected marginal holding and outdating costs can be computed as:
\begin{align*}
H_t(\textbf{x}_t,f_t,q_t)&=\sum\limits_{k=0}^{K-1}h_{t+k}\mathrm{E}\bigg[(q_t-(A_{k,t}+D_{t+k}-\sum\limits_{m=1}^{K-k-1}x_{m,t})^+)^+\bigg|f_t\bigg]\\
&=\sum\limits_{k=0}^{K-1}h_{t+k}\sum\limits_{u=1}^{q_t}\mathrm{P}\bigg((A_{k,t}+D_{t+k}-\sum\limits_{m=1}^{K-k-1}x_{m,t})^+< u\bigg|f_t\bigg)\\
&=\sum\limits_{k=0}^{K-1}h_{t+k}\sum\limits_{u=1}^{q_t}\mathrm{P}\bigg(A_{k,t}+D_{t+k}< u+\sum\limits_{m=1}^{K-k-1}x_{m,t}\bigg|f_t\bigg)\\
&=\sum\limits_{k=0}^{K-1}h_{t+k}\sum\limits_{u=1}^{q_t}R_{k+1}\bigg(u+\sum\limits_{m=1}^{K-k-1}x_{m,t}\bigg),
\end{align*}
and
\begin{align*}
W_t(\textbf{x}_t,f_t,q_t)&=w_{t+K-1}\mathrm{E}[(q_t-A_{K-1,t}-D_{t+K-1})^+|f_t]\\
&=w_{t+K-1}\sum\limits_{u=1}^{q_t}\mathrm{P}(A_{K-1,t}+D_{t+K-1}< u|f_t)\\
&=w_{t+K-1}\sum\limits_{u=1}^{q_t}R_{K,t}(u).
\end{align*}

%%%%%%%%%%%%%%%%%%%%%%%%%%%%%%%%%%%%%%%%%%%%%%%%%%%%%
\subsection{Experiments for the Platelet Inventory Control Problem.}
\label{sec:compound}

We now consider the platelet inventory control problem at Hospital Alpha, as described in $\S$\ref{sec:intro}. At Hospital Alpha, 1) platelets are ordered on a daily basis, and an order placed at the end of the previous day will arrive in the morning of the next day; 2) as demand arises, older products are typically issued first to reduce outdates; and 3) unmet demand is satisfied by emergency deliveries. Therefore, our assumptions for zero lead time, FIFO issuing policy, and lost sales are applicable in this setting.

Platelets have a short lifetime of $K=3$ days, and we consider a planning horizon of 4 weeks (i.e., $T=28$ days). As discussed in $\S$\ref{sec:intro}, we focus on the main source of demand for platelets: cardiac surgeries, and we model daily demand for platelets by a compound Poisson distribution (\citet{gregor1982evaluation, kopach2003models, katsaliaki2008cost}). Similar to two recent studies by \citet{haijema2007blood} and \citet{zhou2011inventory}, we assume that demand over time is independently distributed, but the distribution in different days may not be identical. Based on the cardiac surgery records from January to April in 2014, we identify a significant weekly periodicity, and estimate the average number of surgeries from Monday to Sunday as 2.6, 5.5, 1.9, 3.2, 3.7, 0.1, and 0, respectively. We assume the amount of platelets needed per surgery is stationary; based on the platelet transfusion records, we estimate it as a geometric distribution with mean of 0.32.

Cardiac surgeries are usually scheduled days or even weeks in advance; therefore forecast information on the number surgeries per day is typically available. We consider a forecast horizon equal to product lifetime $K=3$ days, and assume that the forecast is perfect. That is, the number of surgeries at day $t+2$ becomes known at the beginning of day $t$, and will not change at the following days. In this case, at any day $t$, given $f_t$, each of the demand $D_t,...,D_{t+K-1}$ is a sum of a deterministic number of i.i.d. geometric distributions (which is a negative binomial distribution, instead of a compound Poisson distribution). 

Based on the interaction with the blood bank manager at Hospital Alpha, we estimate the unit outdating cost $w$ to be equal to the purchase cost $\$500$. On the other hand, the shortage penalty for blood inventory problems could include the cost of emergent shipment from other blood banks and/or the penalty of postponing the surgeries, which is usually high and often estimated as 2-10 times higher than the purchase cost (\citet{haijema2007blood}). Therefore, we consider three different shortage penalties $p=\$1000,\$2500,\$5000$. Also, we consider a zero holding cost $h=0$ and no discount, i.e., $\beta=1$.

We first benchmark the performances of our policies with the optimal policy solved by dynamic programming. The state of the dynamic program is comprised of a $K-1$ dimensional vector of inventory levels of age $1,...,K-1$, and a $K$ dimensional vector of forecasts on the number of surgeries at days $t,...,t+K-1$. Although the problem size we face here is not too large, it still takes more than 50 hours to compute the optimal policy on a standard 2.6GHz PC, whereas the ordering quantities under our policies can be computed on the fly in an online fashion. On the other hand, while compound Poisson distribution is widely considered in the blood inventory literature (e.g., \citet{gregor1982evaluation, kopach2003models, katsaliaki2008cost}), none of these studies has considered the forecast information on the number of patients per period. A natural question is that how much do we lose by ignoring this information? Therefore, we also compare the performances of our policies with the ``optimal'' policy that does not make use of the forecast information (i.e., it simply treats the demand at each day as a compound Poisson distribution and thus the state of this dynamic program is simply comprised of a $K-1$ dimensional vector of inventory levels of age $1,...,K-1$).

We use $B, TB, OPT$ and $OPT_{wof}$ to denote the marginal-cost dual-balancing policy,\footnote{Since only two cost components $p$ and $w$ are considered here, our marginal-cost dual-balancing policy is the same as both the proportional-balancing policy ($PB$) and dual-balancing policy $(DB)$ proposed in \citet{chao2015approximation}.} the truncated-balancing policy, the optimal policy, and the ``optimal'' policy without forecast information, respectively. For each policy, We generate $10,000$ random scenarios, and use a sample average to estimate the expected total cost. Let $\bar{\mathscr{C}}(\pi)$ and $\bar{\mathscr{C}}(OPT)$ denote the estimated total costs under policies $\pi$ and $OPT$, respectively. We define the performance error of policy $\pi$ as:
$$error(\pi)\mathrel{\mathop:}=\frac{\bar{\mathscr{C}}(\pi)-\bar{\mathscr{C}}(OPT)}{\bar{\mathscr{C}}(OPT)}\times 100\%.$$
We can also characterize the value of forecast information in this setting by assessing the performance improvements of our policies over policy $OPT_{wof}$. Let $\bar{\mathscr{C}}(\pi)$ and $\bar{\mathscr{C}}(OPT_{wof})$ be the estimated total costs under policies $\pi$ and $OPT_{wof}$, respectively. We define the performance improvement of policy $\pi$ as:
$$impr(\pi)\mathrel{\mathop:}=\frac{\bar{\mathscr{C}}(OPT_{wof})-\bar{\mathscr{C}}(\pi)}{\bar{\mathscr{C}}(OPT_{wof})}\times 100\%.$$

\begin{table}[h]
\caption{Performance summary of each policy for the platelet inventory control problem $(w=\$ 500)$}
\vspace{3mm}
\begin{center}
\begin{tabular}{l l r r r r}
\hline
Policy & & \multicolumn{1}{c}{~~~~~B} & \multicolumn{1}{c}{~~~~~TB} & \multicolumn{1}{c}{~~~~OPT} & \multicolumn{1}{c}{~OPT$_{wof}$} \\
\hline
$p=\$ 1000$~~~ & $\bar{\mathscr{C}}$ \hspace{7.0mm} (\$)& ~~~~6813 & ~~~~6684 & ~~~~6174 & ~~7262 \\
& $error$ \hspace{0.2mm} (\%)& 10.4 & 8.3 & 0 & 17.6 \\
& $impr$ \hspace{1.0mm} (\%)& 6.2 & 8.0 & 15.0 & 0 \\
\hline
$p=\$ 2500$ & $\bar{\mathscr{C}}$ \hspace{7.0mm} (\$)& ~~~~10059 & ~~~~9666 & ~~~~8943 & ~~10532 \\
& $error$ \hspace{0.2mm} (\%)& 12.5 & 8.1 & 0 & 17.8 \\
& $impr$ \hspace{1.0mm} (\%)& 4.5 & 8.2 & 15.1 & 0 \\
\hline
$p=\$ 5000$ & $\bar{\mathscr{C}}$ \hspace{7.0mm} (\$)& ~~~~12689 & ~~~~11918 & ~~~~10990 & ~~12999 \\
& $error$ \hspace{0.2mm} (\%)& 15.5 & 8.4 & 0 & 18.3 \\
& $impr$ \hspace{1.0mm} (\%)& 2.4 & 8.3 & 15.5 & 0 \\
\hline
\end{tabular}
\end{center}
\label{table:Table2}
\end{table}

The estimated total cost, performance error, and performance improvement of each policy are reported in Table \ref{table:Table2}. We first observe that both our marginal-cost dual-balancing policy ($B$) and truncated-balancing policy ($TB$) perform significantly better than the theoretical worst-case performance guarantee of two (i.e., error of $100\%$). Further, policy $TB$ has a significant performance improvement over policy $B$, especially when the ratio of unit shortage penalty over unit outdating cost $p/w$ gets large. In particular, we observe from the experiments that when the ratio $p/w$ gets larger, both the optimal ordering quantity and the ordering quantity under policy $B$ gets larger, however the ordering quantity under policy $B$ grows slower than the optimal ordering quantity. In this case, the truncation by lower bound helps correct the under-ordering of policy $B$ and bring the ordering quantity up to a more reasonable level.

Meanwhile, we also observe that the ``optimal'' policy that ignores the forecast information ($OPT_{wof}$) performs poorly, with a performance error of more than 17\% in all three instances, and our policy $TB$ has a substantial performance improvement (more than 8\%) over policy $OPT_{wof}$. Therefore, the value of the forecast information is significant, and implementing an inventory control policy that takes into account such information has a high potential to achieve a better performance in practice.

%%%%%%%%%%%%%%%%%%%%%%%%%%%%%%%%%%%%%%%%%%%%%%%%%%%%%%
\section{Conclusions.}
\label{sec:conclusion}

In this paper, we consider a fixed-lifetime perishable inventory control problem assuming demand is a general stochastic process which can be nonstationary, correlated, and dynamically evolving. Theoretically an optimal ordering policy of this problem can be solved using standard dynamic programming, however it becomes computationally intractable for realistic size problems due to the high dimension of the state space. We first present a computationally efficient algorithm that we call the marginal-cost dual-balancing policy. We then prove that under the marginal-cost accounting scheme, the minimizer of the single-period cost provides a lower bound on the optimal ordering quantity; by combining the specific lower bound we derive and any upper bound on the optimal ordering quantity with the marginal-cost dual-balancing policy, we present a more general class of algorithms that we call the truncated-balancing policy. We prove that when FIFO is an optimal issuing policy, both of our policies have a worst-case performance guarantee of two, i.e., the expected total cost of our policies is at most twice that of an optimal policy. We further provide a necessary and sufficient condition and several easy-to-check sufficient conditions that ensure the optimality of FIFO issuing policy. We also compare our marginal-cost dual-balancing policy with an optimal base-stock policy, and show that the expected total cost of our policy is always at most twice that of an optimal base-stock policy. Finally, we conduct numerical experiments based on a platelet inventory control problem using real data and show that a) our policies perform significantly better than the theoretical worst-case performance guarantee, and b) the truncated-balancing policy significantly outperforms the marginal-cost dual-balancing policy, which illustrates that the lower bound we derive is effective and help improve the performance of the marginal-cost dual-balancing policy.

Our worst-case analysis is built on two novel ideas, the imaginary operation policy and the dynamic unit-matching scheme. In particular, we show that when FIFO is an optimal issuing policy, moving units from older to younger positions in the inventory vector can only decrease the expected total cost. This is very intuitive and helps significantly simplify the analysis by allowing properly modifying the inventory vectors and effectively matching units under two different policies. We believe these ideas are valuable beyond this study and can also be applied to facilitate the analysis for other perishable inventory problems.

% Acknowledgments here
%\ACKNOWLEDGMENT{%
% Enter the text of acknowledgments here
%}% Leave this (end of acknowledgment)
\section*{Acknowledgments.}

The authors thank Prof. Cong Shi and anonymous referees for their valuable suggestions on improving the quality and the presentation of the paper.

% Appendix here
% Options are (1) APPENDIX (with or without general title) or 
% (2) APPENDICES (if it has more than one unrelated sections)
% Outcomment the appropriate case if necessary
%
% \begin{APPENDIX}{<Title of the Appendix>}
% \end{APPENDIX}
%
% or 
%
% \begin{APPENDICES}
% \section{<Title of Section A>}
% \section{<Title of Section B>}
% etc
% \end{APPENDICES}

% References here (outcomment the appropriate case) 

\bibliographystyle{informs2014} % outcomment this and next line in Case 1
\bibliography{paperref}

% CASE 1: BiBTeX used to constantly update the references 
% (while the paper is being written).
%\bibliographystyle{ormsv080} % outcomment this and next line in Case 1
%\bibliography{<your bib file(s)>} % if more than one, comma separated

% CASE 2: BiBTeX used to generate mypaper.bbl (to be further fine tuned)
%\input{mypaper.bbl} % outcomment this line in Case 2

\newpage

\begin{APPENDIX}{~}

\small
%\proof{}
\textbf{Proof of Lemma \ref{lem:c1c2}.}
From the system dynamics, we have $\sum\limits_{k=1}^{K-1}X_{k,t+1}^{\pi}=(Y_t^{\pi}-D_t)^+-(X_{K-1,t}^{\pi}-D_t)^+$, where $(Y_t^{\pi}-D_t)^+$ is the amount of inventory after demand realization at period $t$, and $(X_{K-1,t}^{\pi}-D_t)^+$ is the amount of outdates at period $t$. Then we have:
\begin{align*}
\hat{\mathscr{C}}(\pi)-\mathscr{C}(\pi)=&\sum\limits_{t=1}^{T}\beta^{t-1}\hat{c}\bigg(Q_t^{\pi}+(D_t-Y_t^{\pi})^+-(1-\beta)(Y_t^{\pi}-D_t)^+-\beta (X_{K-1,t}^{\pi}-D_t)^+\bigg)-\beta^{T}\hat{c}\sum\limits_{k=1}^{K-1}X_{k,T+1}^{\pi}\\
=&\sum\limits_{t=1}^{T}\beta^{t-1}\hat{c}\bigg(Q_t^{\pi}+(D_t-Y_t^{\pi})^+-(Y_t^{\pi}-D_t)^+\bigg)+\sum\limits_{t=1}^{T}\beta^{t}\hat{c}\sum\limits_{k=1}^{K-1}X_{k,t+1}^{\pi}-\beta^{T}\hat{c}\sum\limits_{k=1}^{K-1}X_{k,T+1}^{\pi}\\
=&\sum\limits_{t=1}^{T}\beta^{t-1}\hat{c}(D_t-\sum\limits_{k=1}^{K-1}X_{k,t}^{\pi})+\sum\limits_{t=1}^{T-1}\beta^{t}\hat{c}\sum\limits_{k=1}^{K-1}X_{k,t+1}^{\pi}\\
=&\sum\limits_{t=1}^{T}\beta^{t-1}\hat{c}(D_t-\sum\limits_{k=1}^{K-1}X_{k,t}^{\pi})+\sum\limits_{t=2}^{T}\beta^{t-1}\hat{c}\sum\limits_{k=1}^{K-1}X_{k,t}^{\pi}\\
=&\sum\limits_{t=1}^T\beta^{t-1}\hat{c}D_t,
\end{align*}
where the second equality follows from the fact that $\sum\limits_{k=1}^{K-1}X_{k,t+1}^{\pi}=(Y_t^{\pi}-D_t)^+-(X_{K-1,t}^{\pi}-D_t)^+$ as explained above, and the third equality comes from the fact that $(D_t-Y_t^{\pi})^+-(Y_t^{\pi}-D_t)^+=D_t-Y_t^{\pi}$.
\Halmos
%\endproof

\vspace{5mm}

\textbf{Proof of Lemma \ref{lem:im}.}
Without loss of generality, assume that $\tau_i-\tau_1\leq K-1$ (otherwise, we can start from the largest $\tau_i-\tau_j$ that is less than or equal to $K-1$). By construction of policy $IM$, we have $\sum\limits_{k=\tau_i-\tau_1}^{K-1}x_{k,\tau_i}^{IM}=\sum\limits_{k=\tau_i-\tau_1}^{K-1}x_{k,\tau_i}^{B}$ and $x_{k,\tau_i}^{IM}=0, k=\tau_i-\tau_1+1,...,K-1$. Therefore, Inequality \ref{Eq. ine} holds for $k=\tau_i-\tau_1,...,K-1$.

Then, for $j=2,...,i-1$, by construction of policy $IM$, we have $\sum\limits_{k=\tau_i-\tau_{j-1}}^{K-1}x_{k,\tau_i}^{IM}=\sum\limits_{k=\tau_i-\tau_{j-1}}^{K-1}x_{k,\tau_i}^{B}$, $\sum\limits_{k=\tau_i-\tau_j}^{K-1}x_{k,\tau_i}^{IM}=\sum\limits_{k=\tau_i-\tau_j}^{K-1}x_{k,\tau_i}^{B}$ and $x_{k,\tau_i}^{IM}=0, k=\tau_i-\tau_j+1,...,\tau_i-\tau_{j-1}-1$. Therefore, Inequality \ref{Eq. ine} holds for $k=\tau_i-\tau_j,...,\tau_i-\tau_{j-1}-1$.

Finally, by construction of policy $IM$, we have $\sum\limits_{k=\tau_i-\tau_{i-1}}^{K-1}x_{k,\tau_i}^{IM}=\sum\limits_{k=\tau_i-\tau_{i-1}}^{K-1}x_{k,\tau_i}^{B}$ and $x_{k,\tau_i}^{IM}=0, k=1,...,\tau_i-\tau_{i-1}-1$. Therefore, Inequality \ref{Eq. ine} holds for $k=1,...,\tau_i-\tau_{i-1}-1$, which completes the proof.

\vspace{5mm}

%\proof{}
\textbf{Proof of Lemma \ref{lem:ij}.}
We first show that for $t=1,...,T$, $C_{t+1}^{(k)}(\textbf{x}_{t+1},f_{t+1})\leq w/\beta, k=1,...,K-1$, $\forall \textbf{x}_{t+1},f_{t+1}$ such that $\sum\limits_{k=1}^{K-1}x_{k,t+1}< \max\limits_{\tau=1,...,t}\bar{y}_{\tau}-d_{t}$. Suppose for some $t+1$, we have $C_{t+1}^{(k)}(\textbf{x}_{t+1},f_{t+1})> w/\beta$ for some $k=1,...,K-1$ and some $\textbf{x}_{t+1}$ and $f_{t+1}$ such that $\sum\limits_{k=1}^{K-1}x_{k,t+1}< \max\limits_{\tau=1,...,t}\bar{y}_{\tau}-d_t$. At period $t$, let $\textbf{x}_t$ and $q_t$ be such that $x_{k-1,t}=x_{k,t+1}+\epsilon, x_{m-1,t}=x_{m,t+1}, m=1,...,k-1,k+1,...,K-1$, and $x_{K-1,t}=d_t$, where $x_{0,t}=q_t$ and $\epsilon$ is positive but sufficiently small such that $\sum\limits_{k=1}^{K-1}x_{k,t}+q_t=\sum\limits_{k=1}^{K-1}x_{k,t+1}+d_t+\epsilon\leq \max\limits_{\tau=1,...,t}\bar{y}_{\tau}$. Then, FIFO issuing policy will issue $d_t$ units of age $K-1$. Consider another issuing policy $\gamma$ which issues $d_t-\epsilon$ units of age $K-1$ and $\epsilon$ units of age $k-1$. Then, there will be $\epsilon$ more units of outdates under issuing policy $\gamma$ and $\epsilon$ more inventory of age $k$ at the beginning of period $t+1$ under FIFO issuing policy. By assumption, $C_{t+1}^{(k)}(\textbf{x}_{t+1},f_{t+1})> w/\beta$; thus $\gamma$ is strictly better than FIFO, which is a contradiction.

We next show that for $t=1,...,T$, $C_{t+1}^{(k)}(\textbf{x}_{t+1},f_{t+1})\geq 0, k=1,...,K-1, \forall \textbf{x}_{t+1},f_{t+1}$, and $C_{t+1}^{(i)}(\textbf{x}_{t+1},f_{t+1})\leq C_{t+1}^{(j)}(\textbf{x}_{t+1},f_{t+1}), 1\leq i<j\leq K-1$, $\forall \textbf{x}_{t+1},f_{t+1}$ such that $\sum\limits_{k=1}^{K-1}x_{k,t+1}< \max\limits_{\tau=1,...,t}\bar{y}_{\tau}-d_{t}$. The claim is clearly true for $t=T$ since $C_{T+1}(\textbf{x}_{T+1},f_{T+1})=0, \forall \textbf{x}_{T+1},f_{T+1}$. Assume the claim is true for $t+1$. We now show that it is also true for $t$, i.e., $C_{t}^{(k)}(\textbf{x}_{t},f_{t})\geq 0, k=1,...,K-1, \forall \textbf{x}_{t},f_{t}$, and $C_{t}^{(i)}(\textbf{x}_{t},f_{t})\leq C_{t}^{(j)}(\textbf{x}_{t},f_{t}), 1\leq i<j\leq K-1, \forall \textbf{x}_{t},f_{t}$ such that $\sum\limits_{k=1}^{K-1}x_{k,t}< \max\limits_{\tau=1,...,t-1}\bar{y}_{\tau}-d_{t-1}$.

We start with $C_{t}^{(k)}(\textbf{x}_{t},f_{t})\geq 0, k=1,...,K-1, \forall \textbf{x}_{t},f_{t}$. Consider the following two cases. First, suppose at period $t$ we have $\sum\limits_{k=1}^{K-1}x_{k,t}< \max\limits_{\tau=1,...,t}\bar{y}$. Consider the following two systems (both following FIFO issuing policy): System 1 starts from $\textbf{x}_t$ and System 2 starts from $\textbf{x}'_t$, where $x'_{k,t}=x_{k,t}+\epsilon$ and $x'_{m,t}=x_{m,t}, m=1,...,k-1,k+1,...,K-1$, i.e., System 2 has $\epsilon$ more units of age $k$, and $\epsilon$ is positive but sufficiently small such that $\sum\limits_{k=1}^{K-1}x_{k,t}+\epsilon\leq \max\limits_{\tau=1,...,t}\bar{y}_{\tau}$. Let System 2 follow an optimal ordering policy, and let System 1 order $\epsilon$ more units than System 2 and follow an optimal ordering policy afterward. Then, it is sufficient to show that System 1 has no more expected total cost than System 2. Let $\textbf{x}_{t+1}$ and $\textbf{x}'_{t+1}$ be the inventory vectors at period $t+1$ for Systems 1 and 2, respectively. Assume that there are $\xi\leq \epsilon$ more units of outdates in System 2 than in System 1 at period $t$. Then we have $\sum\limits_{k=1}^{K-1}x_{k,t+1}=\sum\limits_{k=1}^{K-1}x'_{k,t+1}+\xi\leq (\max\limits_{\tau=1,...,t}\bar{y}_{\tau}-d_t)^+$, and $\sum\limits_{k=m}^{K-1}x_{k,t+1}\leq\sum\limits_{k=m}^{K-1}x'_{k,t+1}, m=2,...,K-1$. By induction assumption, we have $C_{t+1}^{(i)}(\textbf{x}_{t+1},f_{t+1})\leq C_{t+1}^{(j)}(\textbf{x}_{t+1},f_{t+1})\leq w/\beta, 1\leq i<j\leq K-1,$ $\forall \textbf{x}_{t+1},f_{t+1}$ such that $\sum\limits_{k=1}^{K-1}x_{k,t+1}< \max\limits_{\tau=1,...,t}\bar{y}_{\tau}-d_t$. Therefore, System 1 has no more expected total cost than System 2. Second, suppose at period $t$ we have $\sum\limits_{k=1}^{K-1}x_{k,t}\geq  \max\limits_{\tau=1,...,t}\bar{y}$. Consider the following two systems (both following FIFO issuing policy): System 1 starts from $\textbf{x}_t$ and System 2 starts from $\textbf{x}'_t$, where $x'_{k,t}=x_{k,t}+\epsilon$ and $x'_{m,t}=x_{m,t}, m=1,...,k-1,k+1,...,K-1$, i.e., System 2 has $\epsilon$ more units of age $k$, and $\epsilon$ is any positive number. Let both Systems 1 and 2 follow an optimal ordering policy. Since $\sum\limits_{k=1}^{K-1}x_{k,t}\geq \max\limits_{\tau=1,...,t}\bar{y}$, clearly, the ordering quantities in both systems are zero at period $t$. Let $y_t$ and $y'_t$ be the total inventory levels after ordering in Systems 1 and 2, respectively. Then, $\max\limits_{\tau=1,...,t}\bar{y}\leq y_t\leq y'_t$. Thus the expected cost at period $t$ in System 1 is no more than that in System 2. Let $\textbf{x}_{t+1}$ and $\textbf{x}'_{t+1}$ be the inventory vectors at period $t+1$ for Systems 1 and 2, respectively. Then we have $x_{k,t+1}\leq x'_{k,t+1},k=1,...,K-1$. By induction assumption, we have $C_{t+1}^{(k)}(\textbf{x}_{t+1},f_{t+1})\geq 0,k=1,...,K-1, \forall\textbf{x}_{t+1},f_{t+1}$. Therefore, System 1 has no more expected total cost than System 2

Now it remains to show $C_{t}^{(i)}(\textbf{x}_{t},f_{t})\leq C_{t}^{(j)}(\textbf{x}_{t},f_{t}), 1\leq i<j\leq K-1$, $\forall \textbf{x}_{t},f_{t}$ such that $\sum\limits_{k=1}^{K-1}x_{k,t}< \max\limits_{\tau=1,...,t-1}\bar{y}_{\tau}-d_{t-1}$. Given $\textbf{x}_t,f_t$ such that $\sum\limits_{k=1}^{K-1}x_{k,t}< \max\limits_{\tau=1,...,t-1}\bar{y}_{\tau}-d_{t-1}$, consider the following two systems (both following FIFO issuing policy): System 1 starts from $\textbf{x}'_t$ and System 2 starts from $\textbf{x}''_t$, where $x'_{i,t}=x_{i,t}+\epsilon, x'_{k,t}=x_{k,t}, k\neq i$, and $x''_{j,t}=x_{j,t}+\epsilon, x''_{k,t}=x_{k,t}, k\neq j, 1\leq i<j\leq K-1$, i.e., System 1 starts with $\epsilon$ more units of age $i$ and System 2 with $\epsilon$ more units of age $j$, where $i<j$. Let $\epsilon$ be positive but sufficiently small such that $\sum\limits_{k=1}^{K-1}x_{k,t}+\epsilon\leq \max\limits_{\tau=1,...,t-1}\bar{y}_{\tau}$. Let System 2 follow an optimal ordering policy, and let System 1 order the same amount as System 2 and follow an optimal policy afterward. Then, it is sufficient to show that System 1 has no more expected total cost than System 2. Let $\textbf{x}'_{t+1}$ and $\textbf{x}''_{t+1}$ be the inventory vectors at period $t+1$ in Systems 1 and 2, respectively. Then, we have $x'_{k,t+1}=x''_{k,t+1}, k=1,...,i, x'_{i+1,t+1}\geq x''_{i+1,t+1}$, and $x'_{k,t+1}\leq x''_{k,t+1}, k=i+2,...,K-1$. Since $\bar{y}_t$ minimizes the expected cost at period $t$ and by induction assumption, $C_{t+1}^{(k)}(\textbf{x}_{t+1},f_{t+1})\geq 0, k=1,...,K-1$, $\forall \textbf{x}_{t+1},f_{t+1}$ such that $\sum\limits_{k=1}^{K-1}x_{k,t+1}< \max\limits_{\tau=1,...,t}\bar{y}_{\tau}-d_{t}$, the optimal order-up-to level at each period $t$ is at most $\bar{y}_t$ (because ordering more than $\bar{y}_t$ will increase both the cost at $t$ and the cost-to-go at $t+1$). Assume that there are $\xi\leq \epsilon$ more units of outdates in System 2 than in System 1 at period $t$. Then, $\sum\limits_{k=1}^{K-1}x'_{k,t+1}=\sum\limits_{k=1}^{K-1}x''_{k,t_0+1}+\xi\leq (\max\limits_{\tau=1,...,t}\bar{y}_{\tau}-d_{t})^+$. By induction assumption, we have $C_{t+1}^{(i)}(\textbf{x}_{t+1},f_{t+1})\leq C_{t+1}^{(j)}(\textbf{x}_{t+1},f_{t+1}) \leq w/\beta, 1\leq i<j\leq K-1, \forall\textbf{x}_{t+1},f_{t+1}$ such that $\sum\limits_{k=1}^{K-1}x_{k,t+1}< \max\limits_{\tau=1,...,t}\bar{y}_{\tau}-d_{t}$. Therefore, System 1 has no more expected total cost than System 2, which completes the proof.
\Halmos
%\endproof

\vspace{5mm}

%\proof{}
\textbf{Proof of Lemma \ref{lem:optim}.}
First, since $\bar{y}_t$ minimizes the expected cost at period $t$ and by Lemma \ref{lem:ij}, $C_{t+1}^{(k)}(\textbf{x}_{t+1},f_{t+1})\geq 0, k=1,...,K-1$, $\forall \textbf{x}_{t+1},f_{t+1}$ such that $\sum\limits_{k=1}^{K-1}x_{k,t+1}< \max\limits_{\tau=1,...,t}\bar{y}_{\tau}-d_{t}$, the optimal order-up-to level at each period $t$ is at most $\bar{y}_t$ (because ordering more than $\bar{y}_t$ will increase both the cost at $t$ and the cost-to-go at $t+1$). Therefore, given that we start from zero inventory and an optimal ordering policy is followed at each period under policy $IM$, we have $\sum\limits_{k=1}^{K-1}x_{k,t+1}^{IM}\leq \max\limits_{\tau=1,...,t}\bar{y}_{\tau}-d_{t}, t=1,...,T$. For the case where we start from a high inventory level, by construction, the ordering quantity under policy $IM$ will always be zero until period $t$ such that $\sum\limits_{k=1}^{K-1}x_{k,t}^{IM}\leq \max\limits_{\tau=1,...,t-1}\bar{y}_{\tau}$ (and by forward induction the inequality will continue to hold at all of the following periods), before which no movements of units will be performed since the inventory level under policy $IM$ will be no more than that under policy $B$. Therefore, we always have $\sum\limits_{k=1}^{K-1}x_{k,t+1}^{IM}\leq (\max\limits_{\tau=1,...,t}\bar{y}_{\tau}-d_{t})^+$ at $t+1$ if units are moved at $t$.

We now prove Lemma \ref{lem:optim} in a recursive manner. Recall that for each given sample path, $\mathscr{T}_H=\{\tau_1,...,\tau_n\}$. Consider an variation of policy $IM$, call it $IM_1$; under $IM_1$, the movements of units are only performed at $\tau_1$, and an optimal ordering policy is followed and no movements are performed at the following periods. Then, to show $\mathrm{E}[\mathscr{C}(IM)]\leq \mathrm{E}[\mathscr{C}(OPT)]$, it is sufficient to show $\mathrm{E}[\mathscr{C}(IM_1)]\leq \mathrm{E}[\mathscr{C}(OPT)]$; since if this is true, following a similar argument, the movements at future periods can only further decrease the total cost. Consider any realization of $\tau_1$. Clearly, the total cost under policies $IM_1$ and $OPT$ are the same for all periods $1,...,\tau_1-1$. Without loss of generality, further assume that at $\tau_1$, we have only moved $\epsilon$ units of age $k$ to age zero, $k=1,...,K-1$. Then, after the movements, there are $\epsilon$ more units of age zero but $\epsilon$ fewer units of age $k$ under policy $IM_1$ than under policy $OPT$. 

Consider the following two cases. First, suppose the amount of outdates at $\tau_1$ under policies $IM_1$ and $OPT$ are the same. Then, the total cost at $\tau_1$ under the two policies are the same, and total inventory level at $\tau_1+1$ under the two policies are also the same but the inventory vector under policy $IM_1$ is ``younger'', i.e., $\sum\limits_{k=1}^{K-1}x_{k,\tau_1+1}^{IM}=\sum\limits_{k=1}^{K-1}x_{k,\tau_1+1}^{OPT}$, and $\sum\limits_{k=m}^{K-1}x_{k,\tau_1+1}^{IM}\leq \sum\limits_{k=m}^{K-1}x_{k,\tau_1+1}^{OPT}, m=2,...,K-1$. By Lemma \ref{lem:ij}, we have $C_{t+1}^{(i)}(\textbf{x}_{t+1},f_{t+1})\leq C_{t+1}^{(j)}(\textbf{x}_{t+1},f_{t+1}), 1\leq i<j\leq K-1, \forall\textbf{x}_{t+1},f_{t+1}$ such that $\sum\limits_{k=1}^{K-1}x_{k,t+1}< \max\limits_{\tau=1,...,t}\bar{y}_{\tau}-d_{t}$. Therefore, policy $IM_1$ has no more expected total cost than policy $OPT$.

Second, suppose there are $\xi\leq \epsilon$ more units of outdates at $\tau_1$ under policy $OPT$ than under policy $IM_1$ (this is only possible when we have moved units of age $K-1$ to age zero under policy $IM_1$). Then, we have $\sum\limits_{k=1}^{K-1}x_{k,\tau_1+1}^{IM}=\sum\limits_{k=1}^{K-1}x_{k,\tau_1+1}^{OPT}+\xi$, and $\sum\limits_{k=m}^{K-1}x_{k,\tau_1+1}^{IM}\leq \sum\limits_{k=m}^{K-1}x_{k,\tau_1+1}^{OPT}, m=2,...,K-1$. By Lemma \ref{lem:ij}, we have $C_{t+1}^{(i)}(\textbf{x}_{t+1},f_{t+1})\leq C_{t+1}^{(j)}(\textbf{x}_{t+1},f_{t+1}) \leq w/\beta, 1\leq i<j\leq K-1, \forall\textbf{x}_{t+1},f_{t+1}$ such that $\sum\limits_{k=1}^{K-1}x_{k,t+1}< \max\limits_{\tau=1,...,t}\bar{y}_{\tau}-d_{t}$. Therefore, policy $IM_1$ has no more expected cost than policy $OPT$, which completes the proof.
\Halmos
\endproof

\vspace{5mm}
\textbf{Proof of Theorem \ref{baseb}.}
We prove the theorem in a similar way as for Theorem \ref{optb}, except that now policy $IM$ is constructed based on $BA$ instead of $OPT$. In this case, policy $IM$ also follows a base-stock policy and orders up to the same base-stock level as policy $BA$. 

Recall that for each given sample path, $\mathscr{T}_H=\{\tau_1,...,\tau_n\}$. Consider an variation of policy $IM$, call it $IM_1$; under $IM_1$, the movements of units are only performed at $\tau_1$ and no movements are performed at the following periods. Then, to show $\mathrm{E}[\mathscr{C}(IM)]\leq \mathrm{E}[\mathscr{C}(OPT)]$, it is sufficient to show $\mathrm{E}[\mathscr{C}(IM_1)]\leq \mathrm{E}[\mathscr{C}(OPT)]$; since if this is true, following a similar argument, the movements at future periods can only further decrease the total cost.

Since both policies $IM_1$ and $BA$ follow the same base-stock policy, the total shortage penalty and hold cost under policies $IM_1$ and $BA$ are exactly the same at each period. Consider any realization of $\tau_1$. Clearly, the outdating cost under policies $IM_1$ and $OPT$ are the same for all periods $1,...,\tau_1-1$. Further, since units are only moved from older to younger positions at $\tau_1$ under policy $IM_1$, Lemma \ref{lem:base} implies that with probability one, the total outdating cost under policy $IM_1$ is no more than that under policy $BA$. Therefore, we have $\mathscr{C}(IM_1)\leq \mathscr{C}(BA)$ with probability one. The rest of the proof follows the same way as that for Theorem \ref{optb}.
\Halmos
\vspace{5mm}

%%%%%%%%%%%%%%%%%%%%%%%%%%%%%%%%%%%%%%%%%%%%%%%%%%%%%%%%%%%%%%%%%%%

%\proof{}
\textbf{Proof of Proposition \ref{prop:lower}.}
We start with providing a structural property on the optimal cost-to-go function under the marginal-cost accounting scheme. For $t=1,...,T$, given $\textbf{x}_t$ and $f_t$, let $\tilde{C}_t(\textbf{x}_t,f_t)$ denote the optimal cost-to-go function at period $t$ under the marginal-cost accounting scheme, and as in the paper, let $\Gamma_t(\textbf{x}_t,f_t,q_t)=P_t(\textbf{x}_t,f_t,q_t)+H_t(\textbf{x}_t,f_t,q_t)+W_t(\textbf{x}_t,f_t,q_t)$. Then, the optimality equation under the marginal-cost accounting scheme is:
$$\tilde{C}_t(\textbf{x}_t,f_t)=\min\limits_{q_t\geq 0}\bigg\{\Gamma_t(\textbf{x}_t,f_t,q_t)+\mathrm{E}[\tilde{C}_{t+1}(\textbf{X}_{t+1},F_{t+1})|f_t]\bigg\}.$$
For $k=1,...,K-1$, for the continuous case, let $\tilde{C}_t^{(k)}(\textbf{x}_t,f_t)$ denote the partial derivative of $\tilde{C}_t(\textbf{x}_t,f_t)$ with respect to to $x_{k,t}$; for the discrete case, let $\tilde{C}_t^{(k)}(\textbf{x}_t,f_t)$ denote the incremental of $\tilde{C}_t(\textbf{x}_t,f_t)$ caused by a unit increase of $x_{k,t}$. Then, we have the following result.

\begin{lemma}
\label{lem:lower}
Under Assumption \ref{assump:1}, for $t=1,...,T$, $\tilde{C}_{t+1}^{(k)}(\textbf{x}_{t},f_{t})\leq 0, k=1,...,K-1$, $\forall \textbf{x}_{t},f_{t}$.
\end{lemma}

Proof. 
The claim is clearly true for $t=T$ since $\tilde{C}_{T+1}(\textbf{x}_{T+1},f_{T+1})=0, \forall \textbf{x}_{T+1}, f_{T+1}$. Assume that the claim is true for $t+1$. We now show that it is also true for $t$.

Consider the following two cases. First, suppose at period $t$ we have $\sum\limits_{k=1}^{K-1}x_{k,t}< \max\limits_{\tau=1,...,t}\bar{y}$. Consider the following two systems (both following FIFO issuing policy): System 1 starts from $\textbf{x}_t$ and System 2 starts from $\textbf{x}'_t$, where $x'_{k,t}=x_{k,t}+\epsilon$ and $x'_{m,t}=x_{m,t}, m=1,...,k-1,k+1,...,K-1$, i.e., System 2 has $\epsilon$ more units of age $k$, and $\epsilon$ is positive but sufficiently small such that $\sum\limits_{k=1}^{K-1}x_{k,t}+\epsilon\leq \max\limits_{\tau=1,...,t}\bar{y}_{\tau}$. Let System 1 follow an optimal ordering policy. To define the ordering policy in System 2, let $t_0\in(t,t+K-1]$ be the period such that at all $t,...,t_0-1$, there are still some products that are ordered prior to period $t$ in System 2, while by the beginning of period $t_0$, all of those products are either used to satisfy demand or outdated. Then, we define the ordering policy in System 2 as follows: for each period $t,...,t_0-1$, let System 2 order up to the same level as System 1 (order nothing if this is not feasible), and let System 2 follow an optimal ordering policy afterward.

Then, to prove the lemma, it is sufficient to show that the expected total cost under the marginal-cost accounting scheme in System 2 is no more than that in System 1. By definition of $t_0$, no units ordered at periods $\geq t$ will be outdated by the beginning of period $t_0$. Then, the total cost under the marginal-cost accounting scheme in each system is comprised of the following three parts: i) the shortage penalties that occur at periods $t,...,t_0-1$, ii) the holding costs that occur at periods $t,...,t_0-1$ charged for units ordered at periods $\geq t$, and iii) the total costs (shortage penalties, holding and outdating costs) that occur at periods $\geq t_0$.

i) Consider the shortage penalties that occur at periods $t,...,t_0-1$. By definition of the ordering policy under System 2, after ordering, there is at least the same amount of inventory in System 2 as that in System 1 at each period $t,...,t_0-1$. Therefore, the total shortage penalty at periods $t,...,t_0-1$ in System 2 is no more than that in System 1. 

ii) Consider the holding costs that occur at periods $t,...,t_0-1$ charged for units ordered at periods $\geq t$. Since System 2 started with more inventory, it is possible that for all periods $t,...,t_0-1$, the initial inventory level in System 2 is higher than the total inventory level in System 1 after ordering. Then, the ordering quantity in System 2 is zero for all $t,...,t_0-1$. In this case, System 2 would be empty at the beginning of period $t_0$ and there is nothing to prove. Otherwise, let $s_0\in[t,t_0)$ be the first period such that the ordering quantity in System 2 is strictly positive. Since System 2 started with more inventory than System 1, the amount of outdates in System 2 is at least as much as that in System 1 at each period $t,...,t_0-1$. Therefore, by construction, at all periods $s_0+1,...,t_0-1$, the ordering quantity in System 2 is at least as much as that in System 1, and the total inventory level after ordering in the two systems are the same. Let $\textbf{x}_{t_0}$ and $\textbf{x}'_{t_0}$ be the inventory vectors at period $t_0$ in Systems 1 and 2, respectively. Then, by construction, we have $\sum\limits_{m=k}^{K-1}x'_{m,t_0}\leq \sum\limits_{m=k}^{K-1}x_{m,t_0}, k=1,...,K-1$. Therefore, the holding cost that occurs at periods $t,...,t_0-1$ charged for units ordered at periods $\geq t$ in System 2 is no more than that in System 1.

iii) Consider the total costs that occur at periods $\geq t_0$. At the beginning of period $t_0$, we know that $\sum\limits_{m=k}^{K-1}x'_{m,t_0}\leq \sum\limits_{m=k}^{K-1}x_{m,t_0}, k=1,...,K-1$. By Lemma \ref{lem:ij}, we have $0\leq C_{t+1}^{(i)}(\textbf{x}_{t+1},f_{t+1})\leq C_{t+1}^{(j)}(\textbf{x}_{t+1},f_{t+1}), 1\leq i<j\leq K-1$, $\forall \textbf{x}_{t+1},f_{t+1}$ such that $\sum\limits_{k=1}^{K-1}x_{k,t+1}< \max\limits_{\tau=1,...,t}\bar{y}_{\tau}-d_{t}$. Therefore, the total cost that occurs at periods $\geq t_0$ in System 2 is no more than that in System 1.

Second, suppose at period $t$ we have $\sum\limits_{k=1}^{K-1}x_{k,t}\geq \max\limits_{\tau=1,...,t}\bar{y}$. Consider the following two systems (both following FIFO issuing policy): System 1 starts from $\textbf{x}_t$ and System 2 starts from $\textbf{x}'_t$, where $x'_{k,t}=x_{k,t}+\epsilon$ and $x'_{m,t}=x_{m,t}, m=1,...,k-1,k+1,...,K-1$, i.e., System 2 has $\epsilon$ more units of age $k$, and $\epsilon$ is any positive number. Let System 1 follow an optimal ordering policy. Let both Systems 1 and 2 follow an optimal ordering policy. Since $\sum\limits_{k=1}^{K-1}x_{k,t}\geq \max\limits_{\tau=1,...,t}\bar{y}$, clearly, the ordering quantities in both systems are zero at period $t$. Let $y_t$ and $y'_t$ be the total inventory levels after ordering in Systems 1 and 2, respectively. Then, $\max\limits_{\tau=1,...,t}\bar{y}\leq y_t\leq y'_t$. Thus the expected marginal shortage penalty at period $t$ in System 2 is no more than that in System 1; and there is no marginal holding or outdating cost in either system. Let $\textbf{x}_{t+1}$ and $\textbf{x}'_{t+1}$ be the inventory vectors at period $t+1$ for Systems 1 and 2, respectively. Then we have $x_{k,t+1}\leq x'_{k,t+1},k=1,...,K-1$. By induction assumption, we have $\tilde{C}_{t+1}^{(k)}(\textbf{x}_{t+1},f_{t+1})\leq 0,k=1,...,K-1, \forall\textbf{x}_{t+1},f_{t+1}$. Therefore, under the marginal-cost accounting scheme, System 2 has no more expected total cost than System 1.
\Halmos

With the above result, we now prove the proposition by contradiction. Suppose for some period $t$, given $\textbf{x}_t$ and $f_t$, we have $q_t^{L}> q_t^{OPT}$. Consider a policy $L$, under which $q_t^{L}$ units are ordered at period $t$ and an optimal ordering policy is applied at the following periods. Then, the expected cost-to-go at period $t$ of policy $L$ is $\Gamma_t(\textbf{x}_t,f_t,q_t^{L})+\mathrm{E}[\tilde{C}_{t+1}(\textbf{X}_{t+1}^{L},F_{t+1})|f_t]$, where $\textbf{X}_{t+1}^{L}=\textbf{X}_{t+1}(\textbf{x}_t,q_t^{L},D_t)$. On the other hand, the expected cost-to-go at period $t$ of policy $OPT$ is $\Gamma_t(\textbf{x}_t,f_t,q_t^{OPT})+\mathrm{E}[\tilde{C}_{t+1}(\textbf{X}_{t+1}^{OPT},F_{t+1})|f_t]$, where $\textbf{X}_{t+1}^{OPT}=\textbf{X}_{t+1}(\textbf{x}_t,q_t^{OPT},D_t)$. Since $q_t^{L}> q_t^{OPT}$, by definition of $q_t^{L}$, we have $\Gamma_t(\textbf{x}_t,f_t,q_t^{L})<\Gamma_t(\textbf{x}_t,f_t,q_t^{OPT})$. Further, we have $X_{k,t+1}^{L}\geq_{p} X_{k,t+1}^{OPT}, k=1,...,K-1$ for any realization of $D_{t}$. Therefore, by Lemma \ref{lem:lower}, we have $\tilde{C}_{t+1}(\textbf{X}_{t+1}^{L},F_{t+1})\leq \tilde{C}_{t+1}(\textbf{X}_{t+1}^{OPT},F_{t+1})$ with probability one. Then:
$$\Gamma_t(\textbf{x}_t,f_t,q_t^{L})+\mathrm{E}[\tilde{C}_{t+1}(\textbf{X}_{t+1}^{L},F_{t+1})|f_t]<\Gamma_t(\textbf{x}_t,f_t,q_t^{OPT})+\mathrm{E}[\tilde{C}_{t+1}(\textbf{X}_{t+1}^{OPT},F_{t+1})|f_t],$$
i.e., policy $OPT$ is not optimal for periods $t,...,T$, which is a contradiction.
\Halmos
%\endproof

\vspace{5mm}

%\proof{}
\textbf{Proof of Theorem \ref{opttb}.}
We prove the theorem in a similar way as for Theorem \ref{optb}. The main difference lies in the construction of policy $IM$. In particular, now policy $IM$ is constructed as follows: At each period $t$, given $\textbf{x}_t$ and $f_t$, let the system under policy $IM$ follow an optimal ordering policy. What differentiates policies $IM$ and $OPT$ is that under policy $IM$, at each period after ordering and before demand realization, 1) products in the inventory vector can be ``moved'' from older positions to the position of age 0; and 2) products of age 0 can be intendedly disposed.

At each period $t$, let $y_t^{TB}$ and $y_t^{IM}$ be the total inventory levels after ordering under policies $TB$ and $IM$, respectively. Also, given $\textbf{x}_t^{TB}$ and $f_t$, let $y_t^B$ denote the total inventory level after ordering if the balancing ordering quantity $q_t^B$ is ordered. Then, we partition the set of decision epochs $\{1,...,T\}$ into the following four subsets:
$$\mathscr{T}_P=\{t: y_t^B\geq y_t^{IM}\}, \mathscr{T}_H=\{t: y_t^B< y_t^{IM}, y_t^{TB}= y_t^{B}\},$$
$$\mathscr{T}_{LH}=\{t:y_t^B< y_t^{IM}, y_t^{TB}>y_t^B\}, \mathscr{T}_{UH}=\{t: y_t^B< y_t^{IM}, y_t^{TB}<y_t^B\}.$$

The main objective of constructing policy $IM$ is to bound the the total shortage penalty of policy $TB$ at each $t\in \mathscr{T}_P\cup \mathscr{T}_{UH}$ and the total holding and outdating cost of policy $TB$ charged for the first $q_t^B$ units ordered at each $t\in \mathscr{T}_H\cup \mathscr{T}_{LH}$. In particular, units under policy $IM$ can be moved for $t\in \mathscr{T}_H\cup \mathscr{T}_{LH}\cup \mathscr{T}_{UH}=\{\tau_1,...,\tau_n\}$. The rules of movements are defined in a similar way as before such that after the movements at each $\tau_i$, we have:

(i) There are only positive inventory of age 0 and $\tau_i-\tau_{j}$ under policy $IM$, for all $j=1,...,i-1$ such that $\tau_j\in \mathscr{T}_{H}\cup \mathscr{T}_{LH}$.

(ii) For $j=1,...,i-1$ and $\tau_j\in \mathscr{T}_{H}\cup \mathscr{T}_{LH}$, $\sum\limits_{k=\tau_i-\tau_j}^{K-1}x_{k,\tau_i}^{IM}=x_{\tau_i-\tau_j,\tau_i}^{B}+\sum\limits_{k=\tau_i-\tau_j+1}^{K-1}x_{k,\tau_i}^{TB}$, where $x_{\tau_i-\tau_j,\tau_i}^{B}$ denotes the inventory of age $\tau_i-\tau_j$ at period $\tau_i$ under policy $TB$ if $q_{\tau_j}^B$ instead of $q_{\tau_j}^{TB}$ units are ordered at $\tau_j$.

Note that propoerty (ii) is equivalent to Equation (\ref{eqn:move1}) for $\tau_j\in \mathscr{T}_{H}$ since in that case, we have $q_{\tau_j}^{TB}=q_{\tau_j}^B$.

In addition to movements, we also allow disposals of units at periods in $\mathscr{T}_{UH}$. For $t\in \mathscr{T}_{UH}$, we have $y_{t}^{TB}<y_{t}^B<y_{t}^{IM}$. After the movements of units, there must be at least $y_{t}^{IM}-y_{t}^{TB}$ units of age 0 under policy $IM$. Then, we dispose $y_{t}^{IM}-y_{t}^{TB}$ units of age 0 under policy $IM$ so that after the disposal, we have $y_{t}^{IM}=y_{t}^{TB}$, and none of the above two properties resulted from movements of units is violated. 

Then, similar as before, to show $\mathrm{E}[\mathscr{C}(TB)]\leq 2\mathrm{E}[\mathscr{C}(OPT)]$, it is sufficient to show $\mathrm{E}[\mathscr{C}(IM)]\leq \mathrm{E}[\mathscr{C}(OPT)]$ and $\mathrm{E}[\mathscr{C}(TB)]\leq 2\mathrm{E}[\mathscr{C}(IM)]$, respectively. We have shown in Lemma \ref{lem:optim} that under Assumption \ref{assump:1}, moving units from older to younger positions can only decrease the expected total cost. We now show that disposing units during periods in $\mathscr{T}_{UH}$ can also only decrease the expected total cost. For $t\in \mathscr{T}_{UH}$, since $y_t^{TB}<y_t^B$, by definition of policy $TB$, $y_t^{TB}$ provides an upper bound on the optimal order-up-to level for given $\textbf{x}_t^{TB}$ and $f_t$. Also, similar as before, the inventory vector under policy $IM$ is ``younger'' than that under policy $TB$ after the movements (i.e., for $k=1,...,K-1$, policy $IM$ has no more units of age greater than or equal to $k$). Then it is not difficult to show that the optimal order-up-to level for given $\textbf{x}_t^{IM}$ and $f_t$ is at most $y_t^{TB}$. Therefore, the disposal of inventory from $y_t^{IM}$ to $y_t^{TB}$ will only decrease the expected total cost. Then we have:
\begin{equation}
\label{appendix:imopt}
E[\mathscr{C}(IM)]\leq E[\mathscr{C}(OPT)].
\end{equation}

We next show $\mathrm{E}[\mathscr{C}(TB)]\leq 2\mathrm{E}[\mathscr{C}(IM)]$, which together with Inequality (\ref{appendix:imopt}) lead to our conclusion. By construction of policy $IM$, after the movements and disposals, we have $y_t^B\geq y_t^{IM}, \forall t\in \mathscr{T}_{P}\cup \mathscr{T}_{UH}$. Then clearly:
\begin{equation}
\label{appendix:p}
\sum\limits_{t\in \mathscr{T}_P\cup \mathscr{T}_{UH}}P_t^{B}\leq \sum\limits_{t=1}^TP_t^{IM}.
\end{equation}

Then, define the dynamic unit-matching scheme in a similar way as before, such that the first $q_t^B$ units ordered at each $t\in \mathscr{T}_H\cup \mathscr{T}_{LH}$ under policy $TB$ are matched to units under policy $IM$ on a one to one correspondence, and a matched unit under policy $TB$ stays in inventory no longer than the corresponding unit under policy $IM$. Then, we have:
\begin{equation}
\label{appendix:hlh}
\sum\limits_{t\in \mathscr{T}_H\cup\mathscr{T}_{LH}}H_t^{B}\leq \sum\limits_{t=1}^TH_t^{IM}, \sum\limits_{t\in \mathscr{T}_H\cup\mathscr{T}_{LH}}W_t^{B}\leq \sum\limits_{t=1}^TW_t^{IM}.
\end{equation}

Finally, recall that $\Gamma_t(\textbf{x}_t,f_t,q_t)=P_t(\textbf{x}_t,f_t,q_t)+H_t(\textbf{x}_t,f_t,q_t)+W_t(\textbf{x}_t,f_t,q_t)$. Consider the following three cases. First, suppose $q_t^{TB}=q_t^B$. Then clearly, $\Gamma_t(\textbf{x}_t^{TB},f_t,q_t^{TB})=\Gamma_t(\textbf{x}_t^{TB},f_t,q_t^{B})$. Second, suppose $q_t^{TB}>q_t^B$. Then we have $q_t^{TB}=q_t^L>q_t^B$. Given $\textbf{x}_t$ and $f_t$, it is straightforward to check that $\Gamma_t(\textbf{x}_t,f_t,q_t)$ is convex in $q_t$. Further, since $q_t^{TB}=q_t^L$ minimizes $\Gamma_t(\textbf{x}_t^{TB},f_t,q_t)$, we must have $\Gamma_t(\textbf{x}_t^{TB},f_t,q_t^{TB})\leq \Gamma_t(\textbf{x}_t^{TB},f_t,q_t^{B})$ (This is why the lower bound $q_t^L$ in the definition of policy $TB$ cannot be replaced by tighter ones). Last, suppose $q_t^B>q_t^U$. Then we have $q_t^{TB}=q_t^U$. Since $\Gamma_t(\textbf{x}_t^{TB},f_t,q_t)$ is convex in $q_t$, $q_t^L$ minimizes $\Gamma_t(\textbf{x}_t^{TB},f_t,q_t)$, and $q_t^B>q_t^{TB}\geq q_t^{L}$, we also have $\Gamma_t(\textbf{x}_t^{TB},f_t,q_t^{TB})\leq \Gamma_t(\textbf{x}_t^{TB},f_t,q_t^{B})$. By definition, for any given $f_t$, $\mathrm{E}[P_t^{TB}+H_t^{TB}+W_t^{TB}|f_t]=\Gamma_t(\textbf{x}_t^{TB},f_t,q_t^{TB})$, $\mathrm{E}[P_t^{B}+H_t^{B}+W_t^{B}|f_t]=\Gamma_t(\textbf{x}_t^{TB},f_t,q_t^{B})$. Therefore:
\begin{equation}
\label{appendix:tbb}
\mathrm{E}[P_t^{TB}+H_t^{TB}+W_t^{TB}|f_t]\leq \mathrm{E}[P_t^{B}+H_t^{B}+W_t^{B}|f_t]
\end{equation}

With Inequalities \ref{appendix:p}-\ref{appendix:tbb}, the rest steps to show $\mathrm{E}[\mathscr{C}(TB)]\leq 2\mathrm{E}[\mathscr{C}(IM)]$ are the same as before, which completes the proof. 
\Halmos

%%%%%%%%%%%%%%%%%%%%%%%%%%%%%%%%%%%%%%%%%%%%%%%%%%%%%%%%%%%%%%%%%%%
\vspace{5mm}
%\proof{}
\textbf{Proof of Proposition \ref{prop:fifo}.}
Due to Lemma \ref{lem:ij}, it remains to prove the ``if'' part of the proposition, i.e., if for $t=1,...,T$, $C_{t+1}^{(k)}(\textbf{x}_{t+1},f_{t+1})\leq w/\beta, k=1,...,K-1$, $\forall \textbf{x}_{t+1},f_{t+1}$ such that $\sum\limits_{k=1}^{K-1}x_{k,t+1}< \max\limits_{\tau=1,...,t}\bar{y}_{\tau}-d_{t}$, then Assumption \ref{assump:1} holds. First, since $C_{T+1}^{(k)}(\textbf{x}_{T+1},f_{T+1})=0\leq w/\beta, k=1,...,K-1, \forall \textbf{x}_{T+1},f_{T+1}$, issuing products of age $K-1$ at $T$ clearly results in less cost than issuing younger products and let the oldest products outdate. Further, how we issue products of age less than $K-1$ at $T$ does not affect the total cost. Therefore, Assumption \ref{assump:1} clearly holds for $T$.

Assume that Assumption \ref{assump:1} holds for $t+1$, i.e., starting from period $t+1$, given that $\sum\limits_{k=1}^{K-1}x_{k,t+1}+q_{t+1}\leq \max\limits_{\tau=1,...,t+1}\bar{y}_{\tau}$ at $t+1$ and an optimal ordering policy is implemented at $t+2,...,T$, FIFO is an optimal issuing policy. We now show it also holds for $t$. Starting from period $t$, given $\textbf{x}_t$ and $q_t$ such that $\sum\limits_{k=1}^{K-1}x_{k,t}+q_t\leq \max\limits_{\tau=1,...,t}\bar{y}_{\tau}$, we must have $\sum\limits_{k=1}^{K-1}x_{k,t+1}\leq (\max\limits_{\tau=1,...,t}\bar{y}_{\tau}-d_t)^+$. Thus, under an optimal ordering policy, we have $\sum\limits_{k=1}^{K-1}x_{k,t+1}+q_{t+1}\leq \max\limits_{\tau=1,...,t+1}\bar{y}_{\tau}$. Then, by induction assumption, FIFO is optimal for $t+1,...,T$. It remains to show that FIFO is also optimal at period $t$. Clearly, issuing products of age $K-1$ at period $t$ results in less total cost than issuing younger products and let the oldest products outdate because $C_{t+1}^{(k)}(\textbf{x}_{t+1},f_{t+1})\leq w/\beta, k=1,...,K-1$. Thus, an optimal issuing policy will issue as many oldest products as possible at period $t$. Let $\gamma$ be such an issuing policy. Then, the costs that occur at period $t$ by following FIFO and $\gamma$ are exactly the same. Further, let $\textbf{x}_{t+1}$ and $\textbf{x}'_{t+1}$ be the inventory vectors at period $t+1$ by following FIFO an $\gamma$, respectively. Then, we have $\sum\limits_{k=1}^{K-1}x_{k,t+1}=\sum\limits_{k=1}^{K-1}x'_{k,t+1}$ and $\sum\limits_{k=m}^{K-1}x_{k,t+1}\leq \sum\limits_{k=m}^{K-1}x'_{k,t+1}, m=2,...,K-1$. From the proof of Lemma \ref{lem:ij}, we know that for $t=1,...,T$, $C_{t+1}^{(k)}(\textbf{x}_{t+1},f_{t+1})\leq w/\beta, k=1,...,K-1$, $\forall \textbf{x}_{t+1},f_{t+1}$ such that $\sum\limits_{k=1}^{K-1}x_{k,t+1}< \max\limits_{\tau=1,...,t}\bar{y}_{\tau}-d_{t}$ implies  for $t=1,...,T$, $C_{t+1}^{(i)}(\textbf{x}_{t+1},f_{t+1})\leq C_{t+1}^{(j)}(\textbf{x}_{t+1},f_{t+1}), 1\leq i<j\leq K-1$, $\forall \textbf{x}_{t+1},f_{t+1}$ such that $\sum\limits_{k=1}^{K-1}x_{k,t+1}< \max\limits_{\tau=1,...,t}\bar{y}_{\tau}-d_{t}$. Therefore, FIFO is also optimal at period $t$, which completes the proof.
\Halmos
%\endproof

%%%%%%%%%%%%%%%%%%%%%%%%%%%%%%%%%%%%%%%%%%%%%%%%%%%%%%%%%%%%%%%%%

\vspace{5mm}

%\proof{}
\textbf{Proof of Proposition \ref{prop:demand}.}
Since $\bar{y}_t$ is non-decreasing in $t$, we have $\max\limits_{\tau=1,...,t}\bar{y}_{\tau}=\bar{y}_{t}$. Due to Proposition \ref{prop:fifo}, to show Assumption \ref{assump:1} holds, it is sufficient to show that for $t=1,...,T$, $C_{t+1}^{(k)}(\textbf{x}_{t+1},f_{t+1})\leq w/\beta, k=1,...,K-1$, $\forall \textbf{x}_{t+1},f_{t+1}$ such that $\sum\limits_{k=1}^{K-1}x_{k,t+1}< \bar{y}_{t}-d_{t}$. 

The claim is clearly true for $t=T$ since $C_{T+1}(\textbf{x}_{T+1},f_{T+1})=0, \forall \textbf{x}_{T+1}, f_{T+1}$. Assume that the claim is true for $t+1,...,T+1$. We now show that it is also true for $t$. At period $t$, given $\textbf{x}_t$ and $f_t$ such that $\sum\limits_{k=1}^{K-1}x_{k,t}< \bar{y}_{t-1}-d_{t-1}$, consider the following two systems (both following FIFO issuing policy): System 1 starts from $\textbf{x}_t$ and System 2 starts from $\textbf{x}'_t$, where $x'_{k,t}=x_{k,t}+\epsilon, x'_{m,t}=x_{m,t}, m\neq k, k=1,...,K-1$, i.e., System 2 starts with $\epsilon$ more units of age $k$, and $\epsilon$ is positive but sufficiently small such that $\sum\limits_{k=1}^{K-1}x_{k,t}+\epsilon\leq \bar{y}_{t-1}$. Let System 1 follow an optimal ordering policy, and let System 2 order up to the same level as System 1 at period $t$ (order nothing if this is not feasible) and follow an optimal ordering policy afterward. Then, it is sufficient to show that the expected total cost in System 2 is at most $w\epsilon /\beta$ more than that in System 1.

Let $y_t$ and $y'_t$ be the total inventory levels after ordering in Systems 1 and 2, respectively. Then, by construction, we have $y_t\leq y'_t\leq \bar{y}_t$. Therefore, the total shortage penalty and holding cost at period $t$ in System 2 is no more than that in System 1 (since by definition, $\bar{y}_t$ minimizes the total shortage penalty and holding cost at period $t$ and the sum of shortage penalty and holding cost is convex in ordering quantity). Let $\textbf{x}_{t+1}$ and $\textbf{x}'_{t+1}$ be the inventory vectors at period $t+1$ under Systems 1 and 2, respectively. Then, by construction, we have $x_{1,t+1}\geq x'_{1,t+1}, x_{k,t+1}\leq x'_{k,t+1},k=2,...,K-1$. Assume that there are $\xi\leq \epsilon$ more units of outdates in System 2 than in System 1 at period $t$. Then $\sum\limits_{k=2}^{K-1}x'_{k,t}-\sum\limits_{k=2}^{K-1}x_{k,t}=\epsilon-\xi$. Since $0\leq C_{t+1}^{(k)}(\textbf{x}_{t+1},f_{t+1})\leq w/\beta, k=2,...,K-1$, $\forall \textbf{x}_{t+1},f_{t+1}$ such that $\sum\limits_{k=1}^{K-1}x_{k,t+1}< \bar{y}_{t}-d_{t}$, the expected total cost in System 2 is at most $w\epsilon\leq w\epsilon/\beta$ more than that in System 1, which completes the proof.
\Halmos
%\endproof

\vspace{5mm}

%\proof{}
\textbf{Proof of Proposition \ref{prop:cost}.}
Due to Proposition \ref{prop:fifo}, to show Assumption \ref{assump:1} holds, it is sufficient to show that for $t=1,...,T$, $C_{t+1}^{(k)}(\textbf{x}_{t+1},f_{t+1})\leq w/\beta, k=1,...,K-1$, $\forall \textbf{x}_{t+1},f_{t+1}$ such that $\sum\limits_{k=1}^{K-1}x_{k,t+1}< \max\limits_{\tau=1,...,t}\bar{y}_{\tau}-d_{t}$. 

The claim is clearly true for $t=T$ since $C_{T+1}(\textbf{x}_{T+1},f_{T+1})=0, \forall \textbf{x}_{T+1}, f_{T+1}$. Assume that the claim is true for $t+1,...,T+1$. We now show that it is also true for $t$. At period $t$, given $\textbf{x}_t$ and $f_t$ such that $\sum\limits_{k=1}^{K-1}x_{k,t}< \max\limits_{\tau=1,...,t-1}\bar{y}_{\tau}-d_{t-1}$, consider the following two systems (both following FIFO issuing policy): System 1 starts from $\textbf{x}_t$ and System 2 starts from $\textbf{x}'_t$, where $x'_{k,t}=x_{k,t}+\epsilon, x'_{m,t}=x_{m,t}, m\neq k, k=1,...,K-1$, i.e., System 2 starts with $\epsilon$ more units of age $k$, and $\epsilon$ is positive but sufficiently small such that $\sum\limits_{k=1}^{K-1}x_{k,t}+\epsilon\leq \max\limits_{\tau=1,...,t-1}\bar{y}_{\tau}$. Let System 1 follow an optimal ordering policy, and let System 2 order up to the same level as System 1 at period $t$ (order nothing if this is not feasible) and follow an optimal ordering policy afterward. Then, it is sufficient to show that the expected total cost in System 2 is at most $w\epsilon /\beta$ more than that in System 1.

Let $y_t$ and $y'_t$ be the total inventory levels after ordering in Systems 1 and 2, respectively. Then, by construction, we have $y_t\leq y'_t\leq \max\limits_{\tau=1,...,t}\bar{y}_{\tau}$. Assume that $y'_t-y_t=\eta\leq \epsilon$. Then, there will be at most $\gamma h\eta$ more expected holding cost in System 2 than in System 1 at period $t$. Let $\textbf{x}_{t+1}$ and $\textbf{x}'_{t+1}$ be the inventory vectors at period $t+1$ under Systems 1 and 2, respectively. Then, by construction, we have $x_{1,t+1}\geq x'_{1,t+1}, x_{k,t+1}\leq x'_{k,t+1},k=2,...,K-1$. Assume that there are $\xi\leq \epsilon$ more units of outdates in System 2 than in System 1 at period $t$. Then $\sum\limits_{k=2}^{K-1}x'_{k,t}-\sum\limits_{k=2}^{K-1}x_{k,t}=\epsilon-\xi$. Since $0\leq C_{t+1}^{(k)}(\textbf{x}_{t+1},f_{t+1})\leq w/\beta, k=2,...,K-1$, $\forall \textbf{x}_{t+1},f_{t+1}$ such that $\sum\limits_{k=1}^{K-1}x_{k,t+1}< \bar{y}_{t}-d_{t}$, and $h\leq \frac{1-\beta}{\beta}w$, the expected total cost in System 2 is at most $\gamma h\eta+w\epsilon\leq w\epsilon/\beta$ more than that in System 1, which completes the proof.
\Halmos
%\endproof

\vspace{5mm}

%\proof{}
\textbf{Proof of Proposition \ref{prop:mix}.}
Due to Proposition \ref{prop:fifo}, to show Assumption \ref{assump:1} holds, it is sufficient to show that for $t=1,...,T$, $C_{t+1}^{(k)}(\textbf{x}_{t+1},f_{t+1})\leq w/\beta, k=1,...,K-1$, $\forall\textbf{x}_{t+1},f_{t+1}$ such that $\sum\limits_{k=1}^{K-1}x_{k,t+1}< \max\limits_{\tau=1,...,t}\bar{y}_{\tau}-d_{t}$. 

The claim is clearly true for $t=T$ since $C_{T+1}(\textbf{x}_{T+1},f_{T+1})=0, \forall \textbf{x}_{T+1}, f_{T+1}$. Assume that the claim is true for $t+1,...,T+1$. We now show that it is also true for $t$. At period $t$, given $\textbf{x}_t$ and $f_t$ such that $\sum\limits_{k=1}^{K-1}x_{k,t}< \max\limits_{\tau=1,...,t-1}\bar{y}_{\tau}-d_{t-1}$, consider the following two systems (both following FIFO issuing policy): System 1 starts from $\textbf{x}_t$ and System 2 starts from $\textbf{x}'_t$, where $x'_{k,t}=x_{k,t}+\epsilon, x'_{m,t}=x_{m,t}, m\neq k, k=1,...,K-1$, i.e., System 2 starts with $\epsilon$ more units of age $k$, and $\epsilon$ is positive but sufficiently small such that $\sum\limits_{k=1}^{K-1}x_{k,t}+\epsilon\leq \max\limits_{\tau=1,...,t-1}\bar{y}_{\tau}$. Let System 1 follow an optimal ordering policy, and let System 2 order up to the same level as System 1 at period $t$ (order nothing if this is not feasible) and follow an optimal ordering policy afterward. Then, it is sufficient to show that the expected total cost in System 2 is at most $w\epsilon /\beta$ more than that in System 1.

Let $y_t$ and $y'_t$ be the total inventory levels after ordering in Systems 1 and 2, respectively. Then, by construction, we have $y_t\leq y'_t\leq \max\limits_{\tau=1,...,t}\bar{y}_{\tau}$. Assume that $y'_t-y_t=\eta\leq \epsilon$. The probability that there will be excess inventory after demand realization at period $t$ in either system is upper bounded by $\Phi_t(\max\limits_{\tau=1,...,t}\bar{y}_{\tau})\leq \max\limits_{1<s\leq t\leq T}\Phi_t(\bar{y}_s)=\gamma$. Then, there will be at most $\gamma h\eta$ more expected holding cost and at least $(1-\gamma)p\eta$ less expected shortage penalty in System 2 than in System 1 at period $t$. Let $\textbf{x}_{t+1}$ and $\textbf{x}'_{t+1}$ be the inventory vectors at period $t+1$ under Systems 1 and 2, respectively. Then, by construction, we have $x_{1,t+1}\geq x'_{1,t+1}, x_{k,t+1}\leq x'_{k,t+1},k=2,...,K-1$. Assume that there are $\xi\leq \epsilon$ more units of outdates in System 2 than in System 1 at period $t$. Then $\sum\limits_{k=2}^{K-1}x'_{k,t}-\sum\limits_{k=2}^{K-1}x_{k,t}=\epsilon-\xi$. Since $0\leq C_{t+1}^{(k)}(\textbf{x}_{t+1},f_{t+1})\leq w/\beta, k=2,...,K-1$, $\forall \textbf{x}_{t+1},f_{t+1}$ such that $\sum\limits_{k=1}^{K-1}x_{k,t+1}< \bar{y}_{t}-d_{t}$, and $h\leq \frac{1-\gamma}{\gamma}p+\frac{1-\beta \gamma}{\beta \gamma}w$, the expected total cost in System 2 is at most $\gamma h\eta-(1-\gamma)p\eta+w\epsilon\leq w\epsilon/\beta$ more than that in System 1, which completes the proof.
\Halmos
%\endproof
\end{APPENDIX}

% Appendix here
% Options are (1) APPENDIX (with or without general title) or 
%             (2) APPENDICES (if it has more than one unrelated sections)
% Outcomment the appropriate case if necessary
%
% \begin{APPENDIX}{<Title of the Appendix>}
% \end{APPENDIX}
%
%   or 
%
% \begin{APPENDICES}
% \section{<Title of Section A>}
% \section{<Title of Section B>}
% etc
% \end{APPENDICES}

% Acknowledgments here
%\section*{Acknowledgments.}
% Enter the text of acknowledgments here

% References here (outcomment the appropriate case) 

% CASE 1: BiBTeX used to constantly update the references 
%   (while the paper is being written).
%\bibliographystyle{informs2014} % outcomment this and next line in Case 1
%\bibliography{<your bib file(s)>} % if more than one, comma separated

% CASE 2: BiBTeX used to generate mypaper.bbl (to be further fine tuned)
%\input{mypaper.bbl} % outcomment this line in Case 2

\end{document}